\documentclass[preprint,12pt,3p]{elsarticle}

\usepackage{graphicx}
\usepackage[usenames,dvipsnames]{xcolor}
\usepackage{amsmath}
\usepackage{mathrsfs}
\usepackage{amsfonts}
\usepackage{amssymb}
\usepackage{amsthm}
\usepackage{subcaption}
\usepackage{graphicx}
\usepackage{accents} 
\usepackage{hyperref}

\makeatletter
\newcommand*{\rom}[1]{\expandafter\@slowromancap\romannumeral #1@}
\makeatother

\newcommand\eqs[2]{\eq{#1}-\eq{#2}} 

\newcommand\REMOVE[1]{}

\newcommand\ie{{\it{i.e.~}}}
\newcommand\eg{{\it{e.g.~}}}

\newcommand\rhs{{{\rm r.h.s.~}}}
\newcommand\wrt{{\rm with respect to~}}

\newcommand\zbar{\hat}
\newcommand\pred{{\rm pr}}

\newcommand\exz[1]{}

\newcommand\Lie[1]{{\Lcal}_{#1}\kern-1.65em -\kern 0.55em}
\newcommand\Liex[1]{{\Lcal}_{#1}\kern-1.45em -\kern 0.55em}
\newcommand\dev{\mathrm{dev\,}}

\newcommand\eeb[1]{\eb\left({#1}\right)}

\newcommand\Gms{{\pd_\sigma\Om}}

\def\zerob{\bmi{0}}

\def\ubm{{\bf{u}}}

\def\pbm{{\bf{p}}}

\def\fbm{{\bf{f}}}

\def\gbm{{\bf{g}}}

\def\Kbm{{\bf{K}}}
\def\Mbm{{\bf{M}}}

\def\Abm{{\bf{A}}}
\def\Bbm{{\bf{B}}}
\def\Cbm{{\bf{C}}}

\def\Dbm{{\bf{D}}}

\def\bmi#1{\textbf{\textit{#1}}}

\def\ub{{\bmi{u}}}
\def\vb{{\bmi{v}}}
\def\wb{{\bmi{w}}}
\def\fb{{\bmi{f}}}

\def\nb{{\bmi{n}}}

\def\eb{{\bmi{e}}}

\def\gb{{\bmi{g}}}
\def\hb{{\bmi{h}}}

\def\ab{{\bmi{a}}}
\def\bb{{\bmi{b}}}

\def\rb{{\bmi{r}}}

\def\Ab{{\bmi{A}}}
\def\Bb{{\bmi{B}}}

\def\Fb{{\bmi{F}}}

\def\Sb{{\bmi{S}}}

\def\Kb{{\bmi{K}}}

\def\Ib{{\bmi{I}}}

\def\Lb{{\bf{L}}}
\def\Nbm{{\bf{N}}}

\def\1bm{{\bf{1}}}

\def\Vcal{\mathcal{V}}
\def\Scal{\mathcal{S}}
\def\Zcal{\mathcal{Z}}
\def\Lcal{\mathcal{L}}
\def\Hcal{\mathcal{H}}

\def\vphi{\varphi}
\def\vkappa{\varkappa}

\def\Gopbf{{{\rm C} \kern-0.55em{\bf G}}}%

\def\pd{\partial}
\def\Dop{{{\rm I} \kern-0.2em{\rm D}}}
\def\Aop{{{\rm A} \kern-0.6em{\rm A}}}%
\def\corr{{\rm{cr}}} 

\newcommand{\ol}[1]{\overline{#1}}

\def\RR{{\mathbb{R}}}
\def\Om{\Omega}
\def\what#1{\widehat{#1}}
\def\pdiff#1#2{\frac{\partial {#1}}{\partial {#2}}}
\def\dt#1#2{\frac{{\rm{d}}\,#1}{{\rm{d}}\,#2}}

\def\vphi{\varphi}
\def\lam{\lambda}

\def\dlt{\delta}

\def\phibf{{\mbox{\boldmath$\varphi$\unboldmath}}}
\def\psibf{{\mbox{\boldmath$\psi$\unboldmath}}}

\def\taubf{{\mbox{\boldmath$\tau$\unboldmath}}}
\def\sigmabf{{\mbox{\boldmath$\sigma$\unboldmath}}}

\def\eff{{\rm eff}}
\newcommand{\eq}[1]{(\ref{#1})}
\def\pdt#1{\frac{\partial{#1}}{\partial t}}
\def\dvg{\mbox{\rm div}}

\newcounter{definition}
\newenvironment{mydefinition}[2]%
{\vspace{7pt} \noindent{\bf Definition
  \refstepcounter{definition}\thedefinition. \protect\label{#1}
{\sf({#2})}
\hspace{-2mm}}\rm}%
{} 

\newcounter{remark}
\newenvironment{myremark}[1]%
{\vspace{7pt} \noindent{\bf Remark
  \refstepcounter{remark}\theremark. \protect\label{#1}
\hspace{-2mm}}\rm}%
{\hfill$\bullet$} 

\makeatletter
\def\ps@pprintTitle{%
     \let\@oddhead\@empty
     \let\@evenhead\@empty
     \def\@oddfoot{\footnotesize\itshape
     {\copyright 2017. This manuscript version is made available under the
      \href{https://creativecommons.org/licenses/by-nc-nd/4.0/}{CC-BY-NC-ND 4.0 license}.}
      \hfill}
     \let\@evenfoot\@oddfoot}

\begin{document}

\begin{frontmatter}

\title{Modelling large-deforming fluid-saturated porous media
  using an Eulerian incremental formulation\tnoteref{postprint}}

\author[1]{Eduard Rohan\corref{cor1}}
\ead{rohan@kme.zcu.cz}


\author[1]{Vladim\'{\i}r Luke\v{s}}

\address[1]{European Centre of Excellence, NTIS -- New Technologies for
Information Society Faculty of Applied Sciences, University of West Bohemia,
Univerzitni\'{\i} 22, 30614 Pilsen, Czech Republic}


\cortext[cor1]{Corresponding author.}
\tnotetext[postprint]{The article published in
{\it Advances in Engineering Software}, DOI: \href{http://dx.doi.org/10.1016/j.advengsoft.2016.11.003}{10.1016/j.advengsoft.2016.11.003}.}

\begin{abstract}
The paper deals with modelling fluid saturated porous media subject to large deformation. An Eulerian incremental formulation is derived using the problem imposed in the spatial configuration in terms of the equilibrium equation and the mass conservation. Perturbation of the hyperelastic porous medium is described by the Biot model which involves poroelastic coefficients and the permeability governing the Darcy flow. Using the material derivative with respect to a convection velocity field we obtain the rate formulation which allows for linearization of the residuum function. For a given time discretization with backward finite difference approximation of the time derivatives, two incremental problems are obtained which constitute the predictor and corrector steps of the implicit time-integration scheme. Conforming mixed finite element approximation in space is used. Validation of the numerical model implemented in the SfePy code is reported for an isotropic medium with a hyperelastic solid phase. The proposed linearization scheme is motivated by the two-scale homogenization which will provide the local material poroelastic coefficients involved in the incremental formulation.
\end{abstract}


\end{frontmatter}

\section{Introduction}
Within the small strain theory, behaviour of the
fluid-saturated porous materials (FSPM) was described by M.A.~Biot in \cite{Biot1941}
and further developed in \cite{Biot1955}. It has been proved (see \eg \cite{Auriault1990,Auriault1992} cf. \cite{Rohan-etal-CMAT2015}) that the poroelastic coefficients obtained by
M.A.~Biot and D.G.~Willis in \cite{Biot1957} hold also for quasistatic
situations, cf. \cite{Brown2014}. In the dynamic case the
inertia forces cannot be neglected, however, for small fluid pressure gradients
leading to slow interstitial flows in pores under slow deformation and within the limits of the small strain theory, the effective poroelasticity coefficients of the material describing the macroscopic behaviour can be obtained by homogenization of the static fluid-structure
interaction at the level of micropores, see also \cite{RSW-ComGeo2013}.

Behaviour of the FSPM at finite strains has been studied within the mixture theories \cite{Bowen1976c,Bowen1982} and also using the macroscopic phenomenological approach \cite{Biot1977,Coussy1998}. The former approach has been pursued in works which led to the theory of porous media (TPM) \cite{deBoer1998,Diebels-Ehlers-1996,Ehlers-Bluhm2002}, including elastoplasticity \cite{Borja-CMAME1995}, or extended for multi-physics problems, see \eg \cite{Markert2008}. Micromechanical approaches based on a volume averaging were considered in \cite{deBuhan1998,DormMoliKond-jmps2002}. The issue of compressibility for the hyperelastic solid phase was considered in \cite{Gajo-PRSA2010}, cf. \cite{deBoer1999}.
A very thorough revision of the viscous flows in hyperelastic media with inertia effects was given in \cite{chapelle_moireau_2014}.
There is a vast body of the literature devoted to the numerical modelling of porous media under large deformation, namely in the context of geomechanical applications, thus, taking into account plasticity of the solid phase \cite{carter_randolph_wroth_1979,prevost_1982,nazem_sheng_carter_sloan_2008}. Although the fluid retention in unsaturated media is so far not well understood in the context of deterministic models, a number of publications involve this phenomenon in computational models \cite{meroi_schrefler_zienkiewicz_1995,song_borja_2013}. Concerning the problem formulation, 
three classical approaches can be used. Apart of the total Lagrangian (TL) formulations employed in the works related to the TPM \cite{deBoer1998,Diebels-Ehlers-1996,Ehlers-Bluhm2002}, the updated Lagrangian (UL) formulation has been considered in \cite{meroi_schrefler_zienkiewicz_1995,li_borja_regueiro_2004,shahbodagh-khan_khalili_esgandani_2015}. As an alternative, the ALE formulation provides some advantages, like the prevention of the FE mesh distorsions, cf. \cite{nazem_sheng_carter_sloan_2008}, cf. \cite{Brown2014}.
Transient dynamic responses of the hyperelastic porous medium has been studied in \cite{li_borja_regueiro_2004}, within the framework of the UL formulation with the linearization developed in \cite{Borja-CMAME1995}. In \cite{shahbodagh-khan_khalili_esgandani_2015} the model has been enriched by the hydraulic hysteresis in unsaturated media.

It is worth to recall the generally accepted open problem of ``missing one more equation'' in the phenomenological theories of the FSPM. This deficiency can be circumvented by upscaling. For this, in particular, the homogenization of periodic media can be used, although several other treatments have been proposed, see \eg \cite{Gray-Schrefler-IJNAMG2007} where the Biot's theory was recovered within the linearization concept.
 The conception of homogenization applied to upscale the Biot continuum in the framework of the updated Lagrangian formulation has been proposed \cite{Rohan-etal-2006-CaS}, cf. \cite{Rohan2006}; recently in \cite{Brown2014}, where the ALE formulation was employed and some simplification were suggested to arrive at a computationally tractable two-scale problem.
In \cite{AMC2015-Rohan-Lukes} we described a nonlinear model of the Biot-type poroelasticity  to treat situations when the deformation has a significant influence on the permeability tensor controlling the seepage flow and on the other poroelastic coefficients. 
 Under the small deformation assumption and using the first order gradient theory of the continuum, the homogenized constitutive laws were modified to account for their dependence on the macroscopic variables involved in the upscaled problem. For this, the sensitivity analysis  known from the shape optimization was adopted.

In this paper we propose a consistent Eulerian incremental formulation for the FSPM which is intended as the basis for multi-scale computational analysis of FSPM undergoing large deformation.
In the spatial configuration, we formulate the dynamic equilibrium equation using the concept of the effective stress involved in the Biot model. The fluid redistribution is governed by the Darcy flow model with approximate inertia effects, and by the mass conservation where the pore inflation is controlled through the Biot stress coupling coefficients. The compressibility of the solid and fluid phases are respected by the Biot modulus. The solid skeleton is represented by the neo-Hookean hyperelastic model employed in \cite{li_borja_regueiro_2004}, although any finite strain energy function could be considered.
First, the rate formulation, \ie the Eulerian formulation, is derived by the differentiation
of the residual function \wrt time using the conception of the material derivative consistent with the objective rate principle.
The time discretization leads to the updated Lagrangian formulation whereby the out-of-balance terms are related to the residual function associated with the actual reference configuration.
Although, in this paper,  we do not describe behaviour at the pore level, the proposed formulation
is coherent with the homogenization procedure which can provide the desired effective material parameters at the deformed configuration. In relationships with existing works \cite{Rohan-etal-2006-CaS,Rohan2006} and \cite{Brown2014}, this is an important aspect which justifies  the proposed formulation as an improved modelling framework for a more accurate multi-scale computational analysis  of the porous media.


The paper is organized, as follows.
In Section~\ref{sec:model} we present a model of the hyperelastic FSPM described in the spatial configuration and formulate the nonlinear problem. The inertia effects related to both phases are accounted for in Definition~\ref{def1R}, although in Definition~\ref{def1} we consider the simplified dynamics. Section~\ref{sec:formulation} constitutes the main part of the paper. There we introduce the rate Eulerian formulation using the above mentioned differentiation of the residual function. Then we explain the time discretization which leads to the finite time step incremental formulation suitable for the numerical computation using the finite element (FE) method. In this procedure, we neglect the inertia terms associated with the fluid seepage in the skeleton. Definitions~\ref{def3R} and \ref{def4R} provide two approximate incremental formulations related the predictor and corrector steps, respectively.
In Section~\ref{sec:ex}, we provide a numerical illustration and validation of the incremental formulation using 2D examples considered in the work \cite{li_borja_regueiro_2004}, namely the 1D compression test and the consolidation of a partially loaded layer.
For the validation we use the comparison between responses obtained by the proposed model with those obtained by linear models. Besides inspection of the displacement fields, we show the influence of the time step length on the fluid content, strain energy and the dissipation.

\paragraph{Basic notations}
Through the paper we shall adhere to the following notation.
The spatial position $x$ in the medium is specified through the coordinates
$(x_1,x_2,x_3)$ with respect to a Cartesian reference frame.
The Einstein summation
convention is used which stipulates implicitly that repeated indices are summed
over. For any two vectors $\ab,\bb$, the inner product is $\ab\cdot\bb$, for
    any  two 2nd order tensors $\Ab,\Bb$ the trace of $\Ab\Bb^T$ is $\Ab:\Bb =
    A_{ij}B_{ij}$.
 The gradient \wrt  coordinate $x$  is denoted by $\nabla$.
We shall use gradients of the vector fields; the symmetric gradient of a vector function
    $\ub$, \ie the small strain tensor is $\eeb{\ub} = 1/2[(\nabla\ub)^T + \nabla\ub]$ where the transpose
    operator is denoted by the superscript ${}^T$. The components of a vector $\ub$ will be denoted
by $u_i$ for $i=1,...,3$, thus $\ub = (u_i)$.
To understand the tensorial notation, in this context, $(\nabla \vb)_{ij} = \pd_j v_i$, so that $((\nabla \vb)^T)_{ij} = \pd_i v_j$, where $\pd_i = \pd / \pd x_i$.
By $\ol{D}$ we denote the closure of a bounded domain $D$.
$\nb$ is a unit normal vector defined on a boundary
  $\pd D$, oriented outwards of $D$.
The real number set is denoted by  $\RR$.
Other notations are introduced through the text.

\section{Model of large deforming porous medium}\label{sec:model}
Behaviour of the FSPM is governed by the equilibrium equations and the mass conservation governing the Darcy flow in the large deforming porous material.
Initial configuration associated with material coordinates $X_i$, $i = 1,2,3$ and with domain $\Om_0$ is mapped to the spatial (current) configuration at time $t$ associated with spatial coordinates
$x_i$ and with domain $\Om(t)$. Thus, $x_i = \vphi_i(t,X)$ where $\phibf = (\vphi_i)$ is continuously differentiable. By $\Fb = (F_{ij})$  we denote the deformation gradient given by $F_{ij} = \pd x_i/\pd X_j = \delta_{ij} + \pd u_i/\pd X_j$.

\subsection{Governing equations -- Biot model}\label{sec-biot}

We shall first introduce the governing equations for the problem of fluid
diffusion through hyperelastic porous skeleton.  These involve
the Cauchy stress tensor $\sigmabf = (\sigma_{ij})$, the displacement field $\ub = (u_i)$,  the bulk
pressure\footnote{We adhere to the standard sign convention: the positive pressure means the negative volumetric stress induced by compression.} $p$ and the perfusion velocity  $\wb = (w_i)$ which describes motion of the fluid
relative to the solid phase.  The constitutive laws for $\sigmabf = (\sigma_{ij})$  reads
as
\begin{eqnarray}
\sigmabf &=& -p \Bb + \sigmabf^\eff(\ub)\;, \label{eq:1} \\
\sigmabf^\eff &=& J^{-1}\left[
\mu\bb + (\lam \ln J - \mu)\Ib
\right]\;, \label{eq:1b}
\end{eqnarray}
where $p$ is the pore fluid pressure, $J = \det\Fb$ expresses the relative volume change of the skeleton, $\bb = \Fb\Fb^T$ is the left deformation tensor,
$\Bb = (B_{ij})$ is the Biot coupling coefficient.
Further, by ``eff'' we refer to the effective (hyperelastic)
stress, which is related to a strain energy function. Here we consider a neo-Hookean-type material model employed in \cite{li_borja_regueiro_2004}, which is parameterized by $\mu$ and the Lam\'e constant $\lam$; for a small strain approximation, $\mu' := \mu - \lam \ln J$ expresses the shear modulus.

The inertia of the solid and fluid phases is expressed by the skeleton acceleration $\ab^s$  and by the relative acceleration of the fluid
\begin{equation}\label{eq-aR}
\begin{split}
\ab^R = \dot{\what{\left(\frac{\wb}{\phi}\right)}}\;,\quad \wb = \phi(\vb^f-\vb^s)\;,
\end{split}
\end{equation}
where $\vb^f,\vb^s$ are the solid and fluid phase velocities, respectively\footnote{This expression for $\ab^R$ is merely approximate since the dot means the material derivative \wrt the skeleton, for details see \eg. \cite{chapelle_moireau_2014}.}.
It is now possible to establish the  Darcy law extended for dynamic flows in the moving porous skeleton. According to Coussy \cite{Coussy2004book}
\begin{equation}\label{eq-w}
\begin{split}
\wb & = -\Kb(\nabla p - \rho_f(\fb -  \ab^f -\ab))\\
& \approx -\Kb(\nabla p - \rho_f(\fb - \ab^s - \alpha\ab^R))
\end{split}
\end{equation}
where  $\Kb = (K_{ij})$ is the hydraulic permeability associated with the well known
Darcy law, $\fb= (f_i)$ is the volume force acting on the fluid and $\ab$ is an additional acceleration accounting for the tortuosity effect. The fluid acceleration can be approximated using the one of the skeleton and using the relative acceleration introduced in \eq{eq-aR} multiplied by a correction factor $\alpha$, see also \cite{Coussy2004book}.

Thus, the balance of forces
(equilibrium equation) and the volume conservation read as
\begin{eqnarray}
-\nabla\cdot\sigmabf &=& \ol{\rho}(\gb -\ab^s) - \rho_f \ab^R\;, \label{eq:3R} \\
\Bb:\nabla\pdt \ub + \nabla\cdot\wb + M\pdt p\; &=& 0\;, \label{eq:4R}
\end{eqnarray}
where $M$ is the Biot compressibility, $\gb$ is the volume force acting on the mixture, $\ol{\rho}$ is the mean density.

We remark, that the system \eqs{eq:1}{eq:4R}, first introduced by M. Biot,
see \eg \cite{Coussy2004book,deBoer2000TPM}, is based on the mixture theory. In such approach the
microstructure is represented by volume fractions, whereas any other features
of the microstructural geometry and topology can only be regarded by values of
the constitutive coefficients.  Although, in many realistic situations the
discussed model cannot be recovered by a more rigorous homogenization treatment
(which reflects genuinely the geometry of the microstructure and takes in to account internal friction induced by the fluid viscosity) --- this
phenomenological model can still be used to treat those media for which a detailed description of the microstructure
would be too complex and practically impossible.

\subsection{Nonlinear problem formulation}\label{sec:df}

The differential equations \eqs{eq:3R}{eq:4R} which  hold in $\Om(t)$, are
supplemented by the boundary conditions on $\pd \Om(t)$ which is decomposed into boundary segments according to the following splits:
\begin{equation}\label{eq-cif-2}
\begin{split}
\pd \Om = \pd_u \Om \cup \pd_\sigma \Om\;, \pd_u \Om \cap \pd_\sigma \Om = \emptyset\;,\\
\pd \Om = \pd_p \Om \cup \pd_w \Om\;, \pd_p \Om \cap \pd_w \Om = \emptyset\;.
\end{split}
\end{equation}
In general, we may consider
\begin{equation}\label{eq-cif-3}
\begin{split}
\ub = \ub^\pd \quad \mbox{ on }\pd_u\Om\;,\\
\sigmabf\cdot\nb = \hb \quad \mbox{ on }\pd_\sigma\Om\;,\\
p = p^\pd \quad \mbox{ on }\pd_p\Om\;,\\
\wb\cdot\nb = w_n^\pd \quad \mbox{ on }\pd_w\Om\;.
\end{split}
\end{equation}
 Since the first and the second time derivatives are involved in \eqs{eq-w}{eq:4R}, the initial conditions should involve $\ub,p,\wb$, and also $\dot\ub$.
It is worth noting that these  
imposed in $\Om_0$ must be consistent with the boundary conditions on $\pd\Om_0$ at time $t=0$.. 
On $\pd_w \Om(t)$, the impermeable boundary will be considered for the sake of simplicity, \ie $w_n^\pd = 0$.

In order to establish the weak formulation of the problem defined in Section~\ref{sec:df}, we need the admissibility sets for displacements and for pressures
admissibility sets (we do not specify functional spaces, assuming enough of regularity for all involved unknown functions)
\begin{equation}\label{eq-cif-4}
\begin{split}
V(t) = \{\vb|\; \vb = \ub^\pd \mbox{ on } \pd_u \Om(t)\}\;,\\
W(t) = \{\psibf|\; \nb\cdot\psibf = w_n^\pd \mbox{ on } \pd_w \Om(t)\}\;,\\
Q(t) = \{q|\; q = p^\pd \mbox{ on } \pd_p \Om(t)\}\;,
\end{split}
\end{equation}
and the associated spaces for the test functions
\begin{equation}\label{eq-cif-7}
\begin{split}
V_0(t) = \{\vb|\; \vb = 0 \mbox{ on } \pd_u \Om(t)\}\;,\\
W_0(t) = \{\psibf|\; \nb\cdot\psibf = 0 \mbox{ on } \pd_w \Om(t)\}\;,\\
Q_0(t) = \{q|\; q = 0 \mbox{ on } \pd_p \Om(t)\}\;.
\end{split}
\end{equation}
We consider the force equilibrium equation and the fluid content conservation equation satisfied in the weak sense.
Further we assume that the material properties ensure the following properties:
\begin{equation}\label{eq-material}
\begin{split}
\mu > 0\;,\quad \lam > 0\;,\quad \Kb = \Kb^T\;, \quad \Bb = \Bb^T\;,\quad M\geq 0\;,\\
\end{split}
\end{equation}
whereby $\Kb$ and $\Bb$ are positive definite. 

\begin{myremark}{rem1}
\emph{Incompressibility issue.} When an incompressible solid is considered, the stress decomposition in the solid phase reads $\sigmabf = \sigmabf^\eff - p\Ib$, thus, $\Bb = \Ib$ in \eq{eq:1}. This is the Terzaghi effective stress principle confirmed by both the micromechanical and homogenization approaches, see \eg \cite{deBuhan1998} and \cite{Rohan-etal-CMAT2015}. It is worth noting, that the fluid compressibility has no influence on the stress decomposition. In general, $M$ in \eq{eq:4R} reflects the bulk compressibility of both the phases. For an incompressible solid, $M = \phi/\kappa_f$.
Assuming the solid incompressibility, $J$ expresses the porosity change, thus, the fluid content variation due to the skeleton deformation. Therefore, the strain energy of the hyperelastic solid should be independent of $J$. However, in this paper we employ the same neo-Hookean material used in \cite{li_borja_regueiro_2004} where the strain energy depends on $J$, although in the mass conservation the solid phase is considered as incompressible. We admit this inconsistency to be able to validate the formulation proposed in the present paper using the results obtained in  \cite{li_borja_regueiro_2004}.
It is worth to remark, that the issue of compressible components has been discussed in a detail in \cite{Gajo-PRSA2010}.

\end{myremark}

We can now establish the weak formulation of the dynamic flow-deformation problem which respects the inertia terms induced by the mixture and also by the relative fluid flow in the skeleton, cf. \cite{chapelle_moireau_2014}.

\begin{mydefinition}{def1R}{Weak solution for the three field formulation}
The weak solution of the problem introduced in Section~\ref{sec:df} is the triplet $(\ub,\wb,p) \in V(t)\times W(t)\times Q(t)$ which for all $t > 0$ satisfies
\begin{equation}\label{eq-cif-8R}
\begin{split}
\Phi_t((\ub,\wb,p);(\vb,\psibf,q)) = 0\quad \forall (\vb,\psibf,q) \in V_0(t)\times W_0(t)\times Q_0(t)
\end{split}
\end{equation}
where
\begin{equation}\label{eq-cif-6R}
\begin{split}
\Phi_t((\ub,\wb,p);(\vb,\zerob,0)) & =
\int_{\Om(t)} \sigmabf:\nabla \vb - \int_{\Om(t)} \left( \ol{\rho}(\gb -\ddot\ub) - \rho_f \dot{\what{\phi^{-1}\wb}} \right) \cdot \vb
- \int_\Gms\hb\cdot\vb\;, \\
\Phi_t((\ub,\wb,p);(\zerob,\psibf,0)) & =
\int_{\Om(t)}\psibf\cdot\left(\wb + \Kb\left(\nabla p + \rho_f\left(\ddot\ub + \alpha\dot{\what{\phi^{-1}\wb}} - \fb\right)\right)\right)\;, \\
\Phi_t((\ub,\wb,p);(\zerob,\zerob,q)) & =\int_{\Om(t)}\left (q
\Bb:\nabla\dot \ub - \nabla q  \cdot \wb + M\dot p q
\right ) \;. 
\end{split}
\end{equation}
\end{mydefinition}

The time discretization and the numerical solutions will be considered only for situations when reduced dynamics provides a sufficient physical approximation. By neglecting effects of the accelerations $\ab^R$, the following two field formulation can be obtained.

\begin{mydefinition}{def1}{Weak solution of the two filed formulation --- reduced dynamics}
The weak solution of the problem introduced in Section~\ref{sec:df} is the couple $(\ub,p) \in V(t)\times Q(t)$ which for all $t > 0$ satisfies
\begin{equation}\label{eq-cif-8}
\begin{split}
\Phi_t((\ub,p);(\vb,q)) = 0\quad \forall (\vb,q) \in V_0(t)\times Q_0(t)
\end{split}
\end{equation}
where
\begin{equation}\label{eq-cif-6}
\begin{split}
\Phi_t((\ub,p);(\vb,0)) & =
\int_{\Om(t)} \sigmabf:\nabla \vb - \int_{\Om(t)}  \ol{\rho}( \gb -\ddot\ub)\cdot \vb\quad - \int_\Gms\hb\cdot\vb \;, \\ 
\Phi_t((\ub,p);(\zerob,q)) & =\int_{\Om(t)}\left (q
\Bb:\nabla\dot \ub + \nabla q  \cdot \Kb(\nabla p + \rho_f(\ddot\ub - \fb) + M\dot p q
\right )\;. 
\end{split}
\end{equation}
\end{mydefinition}
In both the Definitions~\ref{def1R} and \ref{def1},
the domain $\Om(t) = \{x = \phibf(t,X), \; X \in \Om_0\}$, such that
$x_i = \vphi_i(t,X) = X_i + u_i(t,x)$.
It reveals the implicit definition of the solution being defined in the domain $\Om(t)$ which itself depends on the solution. To avoid such a complication, the three classical approaches, \ie the TL, UL, or ALE formulations can be used.
 We use the incremental UL (updated Lagrangian) formulation where the recent domain $\Om(t-\dlt t)$ serves as the reference configuration, see \eg \cite{Crisfield2,Simo-Hughes-book1998} for details.

Besides the kinematics associated with the strain filed,
also the solution-dependent volume fraction of the fluid and the fluid density contribute to the nonlinearity of the model equations,
\begin{equation}\label{eq-aux-rho1}
\begin{split}
\rho_f &= \rho_{f0}\exp\{\frac{p - p_0}{\kappa_f}\}\;,\\
\phi^s &=\phi^{s0} J^{-1}\;,\\
\phi \equiv \phi^f &= 1 - \phi^s = 1 - (1-\phi^{f0})J^{-1}\;,\\
\ol{\rho} & = \phi^f \rho_f + \phi^s \rho_s = \phi \rho_f + (1-\phi) \rho_s\;,\\
\ab^R = \dot{\what{\phi^{-1}\wb}} &=
\phi^{-1}(\dot\wb - \wb\phi^{s0}(\phi J)^{-1} \dvg\dot\ub)\;,
\end{split}
\end{equation}
where $p^0$ is the reference pressure and the last expression of $\ab^R$ can be approximated by $\dot\wb / \phi$.
Nevertheless, below we confine our treatment to situations, where $\ab^R$ can be neglected.

The material coefficients $M$, $\Bb$ and $\Kb$ are defined using the following formulas:
\begin{equation}\label{eq-cif-29c-BR}
\begin{split}
M = \frac{\phi}{\kappa_f}\;, \quad \Bb = \Ib\;, \quad \Kb = \ol{\Kb} \exp\{\vkappa(J-1)\}\;,
\end{split}
\end{equation}
where $\kappa_f$ is the bulk modulus of the fluid, $\ol{\Kb}$ is a constant permeability tensor and $\vkappa > 0$.

\section{Consistent incremental formulation}\label{sec:formulation}
In this section we derive the consistent incremental formulation for the model introduced above.
It will be done in two steps.
By time differentiation with a given perturbation velocity field introduced in \eq{eq-cif-11}, associated with the configuration transformation, first we derive the rate of the residual function $\Phi$ which allows to define an approximated problem for computing increments of the displacement and pressure field associated with a time increment $\dlt t$. As a next step, the time discretization is considered. A suitable approximation of the time derivatives occurring in the rate formulation leads to a linear subproblem which can be solved for the increments and, consequently, to updated state variables of the porous medium.

We shall need the perturbation velocity field $\Vcal$ defined in $\Om(t)$ which satisfies
$\Vcal \in V(t)$. This allows us to introduce the perturbed position and the perturbed domain 
\begin{equation}\label{eq-cif-11}
\begin{split}
x' = x + \tau \Vcal(x)\;, \quad \Om' = \Om(t+\tau) = \{x'|\;x \in \Om(t)\}\;,
\end{split}
\end{equation}
where $\tau$ is the perturbation time.
Using the convection velocity $\Vcal$, we can express the associated material derivative $\dot{(~)} = \dt{}{\tau}|_{\tau=0}$ which will be used to derive the rate form of the residuum in \eq{eq-cif-6}.

\begin{myremark}{rem2}
\emph{Material derivative and convection field $\Vcal$.}
The material derivative is introduced by virtue of the perturbation \eq{eq-cif-11}.
It is worth noting that this differentiation is used in the sensitivity analysis of functionals associated with optimal shape problem where the partial differential constraints are considered. In this context, the convection velocity field called ``the design velocity'' is established to link perturbations of the designed boundary with perturbations of the material point in the domain.
The following formulae can be derived easily, see \eg \cite{haslinger_neittaanmaki_1988:_finite_element_approximation,haug_choi_komkov_1986:_design_sensitivity_analysis}.
\begin{equation}\label{eq-B-1}
\begin{split}
\dt{}{\tau}\left(\int_{\Om(t+\tau)} \right)\big |_{\tau = 0} = \int_{\Om(t)}\nabla\cdot \Vcal\;,\\
\dt{}{\tau}\left(\nabla_{x'} v \right)\big |_{\tau = 0} = - \nabla v \nabla \Vcal = -(\pdiff{ v}{x_k} \pdiff{\Vcal_k}{x_i})\;.\\
\end{split}
\end{equation}
\end{myremark}

The incremental formulation which is defined below is based on \eq{eq-cif-8}  rewritten for time $t+\dlt t$, whereby the residual function can be approximated by $\tilde \Phi_{t+\delta t}$, as follows. Let $\zbar t \in [t,t+\delta t]$ and assume the associated configuration (the primary state fields $\Scal(\zbar t) \approx \zbar \Scal= (\zbar \ub,\zbar \wb,\zbar p)$ and all dependent quantities at $t=\zbar t$) can be estimated, while we assume that $\Scal(t)$ is known\footnote{In fact, an approximation $\tilde\Phi_t$ is known rather than an exact residual which should be zero}. In general, denoting by $\Zcal = (\vb,\psibf,q)$ the test functions
\begin{equation}\label{eq-cif-9a}
\begin{split}
\Phi_{t+\delta t}(\Scal(t+\delta t),\Zcal) \approx \tilde \Phi_{t+\delta t}(\tilde\Scal,\Zcal) = \Phi_t(\Scal(t),\Zcal) + \delta\zbar\Phi(\zbar\Scal,\Zcal)\circ (\delta\Scal,\delta t \Vcal)\;,
\end{split}
\end{equation}
where is the increment $\delta\zbar\Phi(\zbar\Scal,\Zcal)\circ (\delta\Scal,\delta t \Vcal)$ is the differential of the residual functional evaluated at time $\zbar t$ \wrt
the state variation $\delta\Scal$ and configuration perturbation $\delta t \Vcal$ associated with the time increment $\dlt t$.
For $\zbar t = t$, \eq{eq-cif-9a} provides a first order approximation of $\tilde \Phi_{t+\delta t}$ using the tangent defined at $t$. In general, a secant-based approximation can be established when an estimate of $\zbar \Scal$ is available at $\zbar t$. This approach is followed in this work; in Section~\ref{sec-t-discr} we introduce a predictor-corrector algorithm within the time discretization scheme applied to \eq{eq-cif-9a}.

\begin{mydefinition}{def2}{Incremental formulation}
Given fields $\Scal(t) = (\ub(t),p(t)) \in V(t)\times Q(t)$ and $\zbar \Scal = (\zbar \ub,\zbar p)$ defined in $\Om(t)$ and $\zbar \Om$, respectively, the new state $\Scal$ at time $t+\delta t$ satisfies
\begin{equation}\label{eq-cif-9}
\begin{split}
 \delta\zbar\Phi(\zbar\Scal;\Zcal)
\circ (\dlt\Scal, \delta t \Vcal) = -\Phi_t(\zbar \Scal;\Zcal)\quad\quad \forall \Zcal \in V_0(t)\times Q_0(t)\;,
\end{split}
\end{equation}
where
\begin{equation}\label{eq-cif-10}
\begin{split}
\ub & = \zbar \ub + \delta \ub\;, p = \zbar p + \delta p\;,\\
\mbox{ where }\quad \delta \ub & = \dot \ub \delta t\;, \quad
\delta p = \dot p \delta t\;,
\end{split}
\end{equation}
pointwise at any $x \in \Om(t)$.
The updated domain $\Om(t+\delta t)$ is obtained using updated coordinates,
$x_i(t+\delta t) = x_i(t) + \dlt u_i(x)$.
\end{mydefinition}

This definition can be applied to introduce a time integration numerical scheme based on the  predictor-corrector steps: after a time discretization, by putting $\zbar t := t$, \eq{eq-cif-9} yields the predictor problem for computing an approximate increment
$\dlt\Scal$ which can be used to introduce an intermediate state $\zbar \Scal$ at $\zbar t \in ]t,t + \dlt t]$.  Then \eq{eq-cif-9} provides a ``corrected approximate increment''.

\subsection{The Lie derivative of the Cauchy stress}
To derive the rate form of the force-equilibrium equation we need the Lie derivative of the Cauchy stress which is related to the 2nd Piola-Kirchhoff stress $\Sb$ by
\begin{equation}\label{eq-cif-5}
\begin{split}
\sigmabf = J^{-1}\Fb \Sb \Fb^T\;,
\end{split}
\end{equation}
recalling the notation $J = \det \Fb$. Since the perturbation of deformation gradient $\Fb$ is $\dlt F_{ij} = \pd \Vcal_i/\pd X_j$, the Lie derivative $\Liex{\Vcal}$ of the Cauchy stress $\sigmabf$
and of the Kirchhoff stress $\taubf = J \sigmabf$ are related to the material derivatives $\dot \sigmabf$ and $\dot\taubf$ by the following expressions
\begin{equation}\label{eq-cif-12}
\begin{split}
\dot\sigmabf & = - \sigmabf \nabla\cdot\Vcal + \nabla\Vcal\sigmabf + \sigmabf(\nabla\Vcal)^T + \Lie{\Vcal} \sigmabf\;,\\
\dot \taubf & = \dot{\what{(J\sigmabf)}} = \nabla\Vcal\taubf + \taubf(\nabla\Vcal)^T + \Lie{\Vcal} \taubf\;,
\end{split}
\end{equation}
where
\begin{equation}\label{eq-cif-12a}
\begin{split}
\Lie{\Vcal} \sigmabf = J^{-1}\Fb\dot\Sb\Fb^T\;,\quad
\Lie{\Vcal} \taubf = \Fb\dot\Sb\Fb^T\;.
\end{split}
\end{equation}
Application of the material derivative \wrt the convection field $\Vcal$, see Remark~\ref{rem2} to differentiate the first integral in \eq{eq-cif-6}$_1$ yields 
(we abbreviate $\Om(t)=\Om$)
\begin{equation}\label{eq-cif-13a}
\begin{split}
& \int_\Om \sigmabf:\nabla \vb \nabla\cdot \Vcal + \int_\Om \dot\sigmabf : \nabla \vb - \int_\Om \sigmabf : (\nabla\vb\nabla\Vcal)
= \int_\Om \left(
\nabla\Vcal \sigmabf + \Lie{\Vcal} \sigmabf \right):\nabla\vb\;,
\end{split}
\end{equation}
where $\Lie{\Vcal} \sigmabf$ is defined using the constitutive relationship.
In the poroelastic medium, the stress is decomposed into its effective and volumetric parts:
$\sigmabf = \sigmabf^\eff - p \Bb$.
The Lie derivative of the effective part of the Cauchy stress can be expressed in terms of the
tangential stiffness tensor and linear velocity strain (\wrt the spatial coordinates, i.e. $e(\dot\ub) = 1/2(\Lb + \Lb^T)$ where $\Lb = \nabla_x \dot\ub$)
\begin{equation}\label{eq-cif-14b}
\begin{split}
\Lie{\Vcal} \sigmabf^\eff = J^{-1}\Fb\dot\Sb\Fb^T = \Dop^\eff \eeb{\dot\ub}\;.
\end{split}
\end{equation}
Above $\Dop^\eff = (D_{ijkl}^\eff)$ is the tangential elasticity tensor which can be identified for a given strain energy function.
The Lie derivative of the volumetric part
$\sigmabf^p = -p \Bb$ can be expressed,
\begin{equation}\label{eq-cif-14}
\begin{split}
\Lie{\Vcal}\sigmabf^p = p (\nabla \Vcal \Bb + \Bb^T(\nabla\Vcal)^T ) - p \Bb \nabla\cdot\Vcal -
\dot{ \what{(p \Bb)}}\;,
\end{split}
\end{equation}
and substituted into \eq{eq-cif-13a} to express the terms involving $\sigmabf^p$, thus
\begin{equation}\label{eq-cif-1-B5}
\begin{split}
& \int_\Om \left(\nabla\Vcal\sigmabf^p + \Lie{\Vcal}\sigmabf^p\right):
\nabla\vb \\
& = - \int_\Om \left(p (\nabla\cdot\Vcal)\Bb:\nabla \vb
+ p \nabla\Vcal\Bb : \nabla\vb - p(\nabla \Vcal \Bb + \Bb^T(\nabla\Vcal)^T )
:\nabla\vb + \dot{ \what{(p \Bb)}} :\nabla\vb
\right)\\
& = - \int_\Om \left(p \Bb:\nabla \vb \nabla\cdot \Vcal
- p(\nabla\Vcal\Bb)^T:\nabla \vb + \dot{ \what{(p \Bb)}}:\nabla\vb
\right)
\;.
\end{split}
\end{equation}

\subsection{Rate formulation in the Euler configuration}
The rate form is introduced to derive a consistent linearization of the nonlinear problem.
We shall need the material derivative of the configuration-dependent variables.
The perturbation defined by \eq{eq-cif-11} induces perturbation of $\phi$, $\rho_f$ and $\ol{\rho}$,
\begin{equation}\label{eq-aux-rho2}
\begin{split}
\dlt\rho_f &= \rho_f\frac{\dlt p}{\kappa_f}\;,\\
\dlt\phi^s &=-\phi^{s0} J^{-1} \dvg\Vcal\;,\\
\dlt\phi  &=  - \dlt\phi^s = \phi^{s0} J^{-1} \dvg\Vcal\;,\\
\dlt\ol{\rho} &= (\rho_f-\rho_s)(1-\phi)\dvg\Vcal + \phi\rho_f\frac{\dlt p}{\kappa_f}\;.
\end{split}
\end{equation}
The following differentials of $\Kb$ and $M$ hold, see \eq{eq-cif-29c-BR}:
\begin{equation}\label{eq-cif-29d-BR}
\begin{split}
\delta M^* & = -\frac{\phi^{s0}}{J\kappa_f}\nabla\cdot\ub^*\;,\\
\delta \Kb^* & = J\vkappa \ol{\Kb} \exp\{\vkappa(J-1)\}\nabla\cdot\ub^*\;.
\end{split}
\end{equation}
Obviously, $\delta \Bb = 0$, since $\Bb = \Ib$ is constant.

In general we consider non-conservative loading by boundary tractions acting on $\Gms$. The classical result of the sensitivity analysis  in shape optimization yields the material derivative of the virtual power induced by the traction forces $\hb$,
\begin{equation}\label{eq-aux-tract}
\begin{split}
\dlt \int_\Gms \hb\cdot\vb =  \int_\Gms (\dlt\hb + (\nb\cdot\Vcal)\Hcal \hb )\cdot\vb\;,
\end{split}
\end{equation}
where $\nb$ is the unit normal vector to $\Gms$ and $\Hcal = k_1 + k_2$ is the 2-multiple of the mean curvature of the deformed surface.

\subsubsection{Incremental operators}

The following expressions
are obtained from \eq{eq-cif-13a}-\eq{eq-aux-rho2}.
We consider a reference state $\Scal = (\ub,\wb, p)$ and its perturbation
$\dlt\Scal = (\dlt\ub,\dlt\wb, \dlt p)$ associated with the convection $\dlt\Vcal$.
Differentiation of $\Phi_t(\Scal;(\vb,\zerob,0))$ yields
\begin{equation}\label{eq-cif-20-BR}
\begin{split}
&\delta\Phi_t(\Scal;(\vb,\zerob,0))
\circ (\dlt\Scal, \dlt\Vcal) = \int_\Om \Dop^\eff e(\dlt \ub):e(\vb) + \int_\Om (\nabla \dlt\Vcal \sigmabf^\eff) :\nabla \vb \\
& - \int_\Om \left(p \Bb:\nabla \vb \nabla\cdot \dlt\Vcal
- p\Bb:\nabla \vb\nabla\dlt\Vcal + (\Bb \dlt p+ \dlt\Bb p):\nabla\vb
\right) \\
& -
\int_\Om \left( \ol{\rho}(\dlt \gb -\dlt \ddot\ub) - \rho_f \dlt( \dot{\what{\phi^{-1}\wb}}) \right) \cdot \vb - \int_\Om \left(\dlt \ol{\rho}(\gb -\ddot\ub) -\dlt \rho_f \dot{\what{\phi^{-1}\wb}}\right) \cdot \vb \\
& - \int_\Om \left( \ol{\rho}(\gb -\ddot\ub) - \rho_f \dot{\what{\phi^{-1}\wb}} \right)  \cdot \vb\nabla\cdot \dlt\Vcal
-\int_\Gms (\dlt\hb + (\nb\cdot\dlt\Vcal)\Hcal \hb )\cdot\vb
\;.
\end{split}
\end{equation}
Differentiation of $\Phi_t(\Scal;(\zerob,\zerob,q))$ yields
\begin{equation}\label{eq-cif-22-BR}
\begin{split}
&\delta\Phi_t(\Scal;(\zerob,\zerob,q))
\circ (\dlt\Scal, \dlt\Vcal) =
\int_\Om q \left(\nabla \cdot \dlt\Vcal\Bb:\nabla\dot \ub + \Bb:\nabla\dlt \dot \ub -
\Bb : \nabla \dot \ub \nabla \dlt\Vcal + \dlt\Bb:\nabla\dot \ub
\right)
\\
& - \int_\Om \left(\wb\cdot \nabla q \nabla \cdot \dlt\Vcal  + \dlt \wb \cdot  \nabla q- \wb \cdot (\nabla q \nabla \dlt\Vcal)\right)  + \int_\Om  \left( q \dlt M \dot p + M \dlt \dot p + q M \dot p \nabla\cdot\dlt\Vcal  \right)
\;.
\end{split}
\end{equation}
In both \eq{eq-cif-20-BR} and \eq{eq-cif-22-BR} we need to express $\wb$ and its perturbation $\dlt\wb$ which can be obtained pointwise in $\Om$,
\begin{equation}\label{eq-cif-21-BR}
\begin{split}
\dlt\wb & = - \Kb\left(\nabla \dlt p -\nabla p \nabla\dlt\Vcal + \rho_f\left(\dlt\ddot\ub + \alpha\dlt(\dot{\what{\phi^{-1}\wb}}) -\dlt \fb\right)
+\dlt\rho_f\left(\ddot\ub + \alpha\dot{\what{\phi^{-1}\wb}} - \fb\right)\right)\\
& \quad -\dlt\Kb \left(\nabla p + \rho_f\left(\ddot\ub + \alpha\dot{\what{\phi^{-1}\wb}} - \fb\right)\right)\;.
\end{split}
\end{equation}


Using \eqs{eq-cif-20-BR}{eq-cif-21-BR} substituted in equation \eq{eq-cif-9}, we obtain the problem for computing increments $(\dlt\ub,\dlt p)$ which are related to a finite time increment $\dlt t$. As the convection field $\dlt\Vcal = \Vcal\dlt t$ is associated with $\dlt\ub$ and since  also the rates $\dot\ub,\ddot\ub$, and $\dot p$ must be expressed in terms of $(\dlt\ub,\dlt p)$,  \eq{eq-cif-9} is not a linear
problem which could be solved directly. To be able to solve this evolutionary problem in time, a suitable time discretization must be introduced, as will be explained in the next section.

\subsection{Time discretization of the incremental subproblem}\label{sec-t-discr}
In the present section we establish a semidiscretized formulation which leads to linear subproblems defined at each time level within the framework of the updated reference configurations and time integration scheme based on the rate formulation treated in the preceding section. Obviously, time discretization of problem \eq{eq-cif-9} is not unique. Implicit Runge-Kutta methods are introduced by means of intermediate time stages, see \eg \cite{park_lee_2009}.
 In computational mechanics, the Newmark-$\alpha$ integration schemes which involve only two consecutive time levels became very popular \cite{li_borja_regueiro_2004,shahbodagh-khan_khalili_esgandani_2015}. When solving a nonlinear problem with dissipation, the naturally arising but seldom considered question concerns with the influence of the time step length and the number of iterations employed to solve the nonlinear problem by the Newton-Raphson method on the dissipated energy.


Here we consider a multistep approach based on combination of forward and backward finite differences to discretize the time derivatives and the convection $\dlt t \Vcal$; optionally a two-level computational scheme of the ``leap frog'' type, \cite{park_lee_2009}
can be defined which uses intermediate time levels.
For this we introduce a uniform time discretization of the time interval $[0,T]$ upon introducing the
time levels $t_k$, $k = 0,1,\dots$, such that $\bigcup_{k=1}^{\bar k}[t_{k-1},t_k] = [0,T]$. By $a^k$ we refer to an approximated value of $a(t_k)$. Further, we define $\delta \ub^k=\ub^k - \ub^{k-1}$ and consider the following approximation of time derivatives at $t_* \in [t_{k-1},t_k]$,
\begin{equation}\label{eq-cif-31aR}
\begin{split}
\dot \ub(t_*) & \approx (\ub^k - \ub^{k-1})/\delta t = \delta \ub^k / \delta t\;,\\
\ddot \ub(t_*) & \approx (\delta \ub^k -  \delta \ub^{k-1})/(\delta t )^2 =
(\ub^k - 2\ub^{k-1} + \ub^{k-2})/(\delta t )^2 \;\\
\dddot\ub(t_*) & \approx 
(\delta \ub^k - 2\delta \ub^{k-1} + \delta \ub^{k-2})/(\delta t )^3 \\
        & = (\ub^k - 3\ub^{k-1} + 3\ub^{k-2} - \ub^{k-3})/(\delta t )^3\;.
\end{split}
\end{equation}

Below we propose a predictor-corrector time-integration scheme, whereby the corrector uses the rate of the residual at an interpolated time level.
The need for the 3rd time derivative in the rate form equations (in fact $\delta \ddot\ub$ is involved) apparently induces the 4 step numerical scheme, see \eq{eq-cif-31aR}$_4$, however, as will be shown later, the predictor step circumvents this complication so that 3 consecutive steps are involved only.
In order to get a linear subproblem arising from \eqs{eq-cif-20-BR}{eq-cif-21-BR}, the convection velocity field $\Vcal$ involved in these expressions can be approximated by the backward difference, \ie $\Vcal(t_k) \approx (\ub^k - \ub^{k-1})/\delta t$, while the time derivatives are approximated using \eq{eq-cif-31aR} with $t_* = t_k$. For the corrector step, $\Vcal$ will be approximated using an estimate of $(\ub^{k+1} - \ub^{k})/\delta t$ computed by the predictor.

In the rest of this paper, we shall confine to the simplified dynamics of the porous medium which disregards the  inertia induced by the relative fluid-solid motion. The time discretization of \eq{eq-cif-6} yields the following approximation
\begin{equation}\label{eq-cif-23-BR}
\begin{split}
 \tilde\Phi_{k}(\{\Scal^l\}_{l\leq k};\Zcal) & =
 \int_{\Om_k}\left (
\sigmabf^k: \nabla \vb - \ol{\rho}^k(\gb^k -(\delta\ub^k-\delta\ub^{k-1})/(\delta t )^2)\cdot \vb
\right) \\
& \quad + \int_{\Om_k}\left (\Bb:\nabla\frac{\delta\ub^k}{\delta t} q
- \wb^k\cdot \nabla q + q M^k \frac{\delta p^k}{\delta t}\right) 
\;,\\
 \wb^k & = -\Kb^k(\nabla p^k- \rho_f^k(\fb^k - (\delta\ub^k-\delta\ub^{k-1})/(\delta t )^2))
\end{split}
\end{equation}
for all $\Zcal = (\vb,q) \in V_0(t_k)\times Q_0(t_k)$.

The discretization of the incremental operators \eqs{eq-cif-20-BR}{eq-cif-21-BR} evaluated at time $\zbar t \in [t_k,t_{k+1}]$ with $\dlt\Vcal = \dlt\ub^*$ yields expressions where we use the following notation
$\Om(\zbar t) =\zbar \Om, \zbar \Dop^\eff := \Dop^\eff(\zbar t), \zbar \sigmabf^\eff := \sigmabf^\eff(\zbar t)$ and $\zbar p = p(\zbar t)$; all other coefficients are evaluated according to the configuration $\zbar \Om$.
From \eq{eq-cif-20-BR} we get
\begin{equation}\label{eq-cif-20aR}
\begin{split}
&\delta\tilde\Phi_k(\zbar\Scal,\{\Scal^l\}_{l\leq k+1};(\vb,0))
\circ (\delta \ub, \delta p,\delta \gb) = \\
&\int_{\zbar\Om} \zbar\Dop^\eff e(\delta\ub^{k+1}):e(\vb) + \int_{\zbar\Om} (\nabla \delta \ub^{k+1} \zbar\sigmabf^\eff(t_k)) :\nabla \vb \\
& - \int_{\zbar\Om} \zbar p \left[\Bb:\nabla\vb(\nabla\cdot \delta\ub^*)
- \nabla\vb \nabla\delta\ub^* : \Bb\right] - \int_{\zbar\Om}(\delta p^{k+1} \Bb:\nabla\vb
+ \zbar p\dlt\Bb:\nabla\vb )\\
& -
\int_{\zbar\Om}  \zbar{\ol{\rho}}\left(\dlt \gb^{k+1} -(\delta \ub^{k+1} - 2\delta \ub^{k} + \delta \ub^{k-1})/(\dlt t)^2 \right)\cdot\vb\\
& -\int_{\zbar\Om} \left( \gb^k -(\delta \ub^{k+1} - \delta \ub^{k})/(\dlt t)^2 )\right) \zbar R(\dlt\ub^*,\dlt p^*) \cdot \vb
-\int_{\zbar{\Gms}} (\dlt\hb^{k+1} + (\nb\cdot\delta \ub^{k})\Hcal \hb )\cdot\vb
\;.
\end{split}
\end{equation}
Eq. \eq{eq-cif-21-BR} is multiplied by $\delta t $, so that  $\dot\ub$ in \eq{eq-cif-6}$_2$ is replaced by $\delta \ub$.
The discretization results in
\begin{equation}\label{eq-cif-21aR}
\begin{split}
& \delta\tilde\Phi_k(\zbar\Scal,\{\Scal^l\}_{l\leq k+1};(0,q))
\circ (\delta \ub, \delta p,\delta \fb)\delta t  = \\
& \int_{\zbar\Om} \left(q (\Bb:\nabla\delta \ub^{k+1})(\nabla\cdot \delta \ub^*)
- q \nabla \delta \ub^{k+1} \nabla \delta \ub^* : \Bb + q \Bb:\nabla(\delta \ub^{k+1} - \delta \ub^k) + q\delta\Bb \nabla\dlt\ub^*
\right)\\
& - \delta t \int_{\zbar\Om} \left(\delta\wb^{k+1} \cdot\nabla q + {\wb^{k+1}\cdot [\nabla q \nabla \cdot \delta \ub^* - \nabla q \nabla \delta \ub^*]}\right)
\\
& +\int_{\zbar\Om}q\left(M(\dlt p^{k+1}-\dlt p^k) + \dlt M  \dlt p^k + M \dlt p^{k+1}(\nabla\cdot\dlt \ub^*) \right)
\;,
\end{split}
\end{equation}
where the differential of the seepage velocity is given by the following approximation:
\begin{equation}\label{eq-cif-21bR}
\begin{split}
\delta\wb^{k+1} = & - \delta\Kb(\nabla p^{k} - \rho_f^k(\fb^k - (\delta \ub^{k+1} -\delta \ub^{k})/(\dlt t)^2))\\
&+\Kb^k (\nabla p^{k} \nabla \delta \ub^{k} +\frac{\rho_f^k}{\kappa_f}\delta p^{k}(\fb^k - (\delta \ub^{k+1} -\delta \ub^{k})/(\dlt t)^2))\\
&-\Kb^k\left(\nabla \delta p^{k+1}-\rho_f^k\left(\delta \fb^{k+1} - (\delta \ub^{k+1} - 2\delta \ub^{k} + \delta \ub^{k-1})/(\dlt t)^2\right)\right)\;.
\end{split}
\end{equation}
Above $\delta\ub^*$ can be substituted by $\delta\ub^k$.

In order to simplify notation, we shall introduce the following tensors and material coefficients depending on the solution at time $\zbar t$ and on other variables associated with time~$\zbar t$:
\begin{equation}\label{eq-cif-26-BR}
\begin{split}
\tilde\Aop &= (\tilde A_{ijkl})\;,\quad \tilde A_{ijkl} = D_{ijkl}^\eff + \sigma_{lj}^\eff\delta_{ki} + \hat p {(B_{il}\delta_{jk} - B_{ij}\delta_{kl})} \;,\\
\tilde\Bb(\vb) &= \Bb (\nabla\cdot\vb - (\nabla\vb)^T)\;,\\
\tilde\Kb(\vb) &= (\nabla\cdot\vb) \Kb - \Kb(\nabla\vb)^T - (\nabla\vb)\Kb\;,\\
\what\Kb(\vb,q) & = \left[\Ib({\kappa_f}^{-1}q + \nabla\cdot\vb) - \nabla\vb\right]\Kb\,\\
\tilde M(\vb) &= M \nabla\cdot\vb\;,\\
\tilde R(\vb,q) & = \left(\ol{\rho}
+ (\rho_f-\rho_s)(1-\phi)\right)\nabla\cdot\vb + \frac{\phi\rho_f}{\kappa_f}q  = \rho_f\left(\nabla\cdot\vb + \frac{\phi}{\kappa_f}\right) q\;.
\end{split}
\end{equation}
Further, we shall employ an abbreviated notation:
\begin{equation}\label{eq-cif-32}
\begin{split}
\delta \ub^{k+1} \mapsto \ub\;,\quad \delta p^{k+1}\mapsto p\;, \\
\delta \ub^k \mapsto \bar\ub\;,\quad \delta p^{k}\mapsto \bar p\;.
\end{split}
\end{equation}


Two incremental linearized subproblems will be formulated which
provide the solutions of the problem introduced in the Definition~\ref{def2}. In this context, we emphasize the different meaning of $\ub$ and $p$ which now denote the increments $\dlt \ub$ and $\dlt p$.
By $\delta V$ and $\delta Q$ we denote sets of admissible increments $\delta \ub^{k+1}\equiv\ub$ and $\delta p^{k+1}\equiv p$ respectively.
For a given time discretization with a time step $\dlt t$,
the sets $\dlt V(t_{k+1}) \times \dlt Q(t_{k+1})$ and $V_0(t_k)\times Q_0(t_k)$ are defined according to the Dirichlet boundary conditions considered according to \eq{eq-cif-31aR}. Also the tensors $\Aop,\Bb,\Kb$ and $M$ are defined in $\Om \equiv \Om(t_k)$, and the volume load $\gb^{k+1}$ is given at~$t_k$.

\subsubsection{The predictor step}\label{sec-pred}
The predictor step is based on the approximation \eq{eq-cif-9a} with $\dlt \zbar \Phi$ evaluated at $t_k$
\begin{equation}\label{eq-pred1}
\begin{split}
0 \stackrel{!}{=} \Phi_{k+1}(\Scal(t_{k+1}),\Zcal) \approx \Phi_k(\{\Scal^l\}_{l\leq k},\Zcal) + \delta\Phi_k(\{\Scal^l\}_{l\leq k},\Zcal)\circ (\delta\Scal,\delta \Vcal)\;,
\end{split}
\end{equation}
where  $\delta \Vcal = \bar\ub$ is approximated by $\dlt \ub^k$, thus, recalling \eq{eq-cif-32}.
For all quantities associated with the configuration at $t_k$ we use the following simplified notation:
$\Om := \Om(t_k), \Dop^\eff := \Dop^\eff(t_k), \sigmabf^\eff := \sigmabf^\eff(t_k)$; all other coefficients are evaluated according to the configuration $\Om_k= \Om(t_k)$.

By virtue of Definition~\ref{def2}, the following is to be solved: Find $(\ub,p) \in \delta V(t_{k+1}) \times \delta Q(t_{k+1})$,  such that
\begin{equation}\label{eq-pred2}
\begin{split}
 &\int_{\Om} \Dop^\eff e(\ub):e(\vb) + \int_{\Om} \sigmabf^\eff :\nabla \vb(\nabla \ub)^T \\
& - \int_{\Om}  p^k \Bb:\left[(\nabla\vb)(\nabla\cdot \ub)
- \nabla\vb \nabla\ub \right]
- \int_{\Om} (p \Bb +  p^k \delta \Bb):\nabla\vb \\
& -\int_\Om \left(\gb^k -(\ub - \bar\ub)/(\dlt t)^2 )\right)\tilde R(\ub^*,p^*) \cdot \vb
+\int_\Om  \ol{\rho}^k\ub(\dlt t)^{-2} \cdot \vb \\
& =
\int_\Om  \ol{\rho}^k\left(\dlt \gb^{k+1} +  (2\bar \ub - \bar\ub^{k-1})/(\dlt t)^2 \right)\cdot\vb
+ \int_\Gms (\dlt\hb^{k+1} + (\nb\cdot\ub^*)\Hcal \hb^k )\cdot\vb)\\
& +  \int_{\Gms} \hb^k\cdot\vb -\int_{\Om}\left (
\sigmabf^k: \nabla \vb - \ol{\rho}^k(\gb^k -(\bar\ub-\bar\ub^{k-1})/(\delta t )^2)\cdot \vb
\right)\;,
\end{split}
\end{equation}
for all $\vb \in V_0(t_k)$ and 
\begin{equation}\label{eq-pred3}
\begin{split}
& \int_{\Om} q \Bb:\left[(\nabla\ub)(\nabla\cdot \ub^*)
- \nabla\ub \nabla\ub^* \right] + \int_{\Om} q \Bb:\nabla\ub
+ \int_{\Om} q M\left( 1 + \nabla\cdot \ub^*\right) p + \int_{\Om} q \delta M^* p
\\
&
+ \delta t\int_{\Om}  \Kb\left(\nabla p +\rho_f^k (\ub - 2\bar\ub + \bar\ub^{k-1})/(\dlt t)^2\right) \cdot \nabla q
\\
& -  \delta t\int_{\Om} \Kb \left(\nabla p^k\nabla \ub^* +\frac{\rho_f^k}{\kappa_f} p^*(\fb^k - (\ub -\bar\ub)/(\dlt t)^2))\right) \cdot \nabla q \\
&+ \delta t\int_{\Om} \Kb \left(\nabla p^k -\rho_f^k(\fb^{k} -(\ub - \bar\ub)/(\delta t)^2  )
\right)\cdot (\nabla q\nabla \cdot \ub^*- \nabla q \nabla \ub^*)\\
&+ \delta t\int_{\Om} \delta \Kb^* (\nabla p^k -\rho_f^k(\fb^{k} -(\ub - \bar\ub)/(\delta t)^2  ))\cdot \nabla q \\
& =  \int_{\Om} q\left(\Bb:\nabla \bar\ub - \delta\Bb^*:\nabla\ub^*\right) + \int_{\Om} q M \bar p \\
& \quad
-\int_{\Om}\left (q\Bb:\nabla{\bar\ub}
+ \delta t \Kb(\nabla p^k -\rho_f^k({\fb^{k+1}} -(\bar\ub - \bar\ub^{k-1})/(\delta t)^2  ))\cdot \nabla q  + q M \bar p\right)
\;,
\end{split}
\end{equation}
for all $q \in Q_0(t_k)$.

\begin{mydefinition}{def3R}{Weak solution of the predictor subproblem}
The weak solution of the predictor subproblem at time level $t_{k+1}$ is the couple
$(\ub,p) \in \delta V(t_{k+1}) \times \delta Q(t_{k+1})$ which satisfies
\begin{equation}\label{eq-pred6}
\begin{split}
 &
\int_{\Om}\left(\ol{\rho}^k + \tilde R(\ub^*,p^*)\right)/(\dlt t)^{2}\ub \cdot \vb
+\int_{\Om}[\tilde \Aop \nabla \ub]:\nabla\vb  - \int_{\Om} p \Bb:\nabla\vb\\
&  =
\int_\Om \left( \ol{\rho}^k\ + \tilde R(\ub^*,p^*)\right)(\bar \ub /(\dlt t)^2 + \gb^{k+1})\cdot\vb  \\
&+\int_\Gms (\hb^{k+1} + (\nb\cdot\ub^*)\Hcal \hb^k )\cdot\vb
-\int_{\Om}\left (
\sigmabf^k - p^k \delta\Bb  \right):\nabla\vb \;,
\end{split}
\end{equation}
for all $\vb \in V_0(t_k)$ and
\begin{equation}\label{eq-pred7}
\begin{split}
& \int_{\Om} q [\Bb + \tilde\Bb(\ub^*)]:\nabla\ub
 +  \delta t\int_{\Om} \Kb \nabla p \cdot \nabla q
 + \int_{\Om} q (M + \tilde M(\ub^*) + \delta M^*) p\\
&+  (\dlt t)^{-1}\int_{\Om} \rho_f^k \nabla q \cdot\left(\Kb + \what\Kb(\ub^*,p^*) + \dlt\Kb\right)\ub
\\
& =   -\int_{\Om} q\left(\delta\Bb :\nabla \bar\ub \right)  
- \delta t \int_{\Om}
\left (\Kb + \tilde\Kb(\ub^*) +  \delta \Kb\right)\nabla p^k \cdot \nabla q\\
&
+ \delta t \int_{\Om}\rho_f^k\nabla q \cdot\left(\Kb + \what\Kb(\ub^*,p^*) + \delta \Kb\right)
(\fb^{k} + \bar\ub /(\delta t)^2  )
\;,
\end{split}
\end{equation}
for all $q \in Q_0(t_k)$,
where $\ub^* = \bar\ub$. 

\end{mydefinition}

It is worth noting that the final form of the predictor subproblem involves the backward differences $\bar\ub$ only, since
$\bar\ub^{k-1}$ and $\bar p$ can be eliminated form \eq{eq-pred2} and \eq{eq-pred3} due to the consistent approximation of $\Phi_{k+1}$ using $\Phi_k$ and the differential $\dlt\Phi_k$, both related to the same configuration $\Om_k$.




\subsubsection{The corrector step}
The predictor step provides the increment $(\dlt\ub^{k+1}, \dlt p^{k+1})$ approximation  denoted as $(\ub^\pred,p^\pred)$ which is related to the time step $\dlt t$. It enables
to establish the configuration $\zbar \Om$ at time $\zbar t = t_{k+1/2} = (t_k + t_{k+1})/2$ and to define the secant approximation of the residual functional $\Phi_{k+1}(\Scal(t_{k+1}),\Zcal)$ in the sense of Definition~\ref{def2}. The displacement and pressure fields are introduced, as follows: $\zbar \ub = \ub^k + \frac{1}{2}\ub^\pred$ and $\zbar p = p^k + \frac{1}{2}p^\pred$.
For all quantities associated with the updated configuration at $\zbar t$ we use the following simplified notation $\Om := \Om(\zbar t) =\zbar \Om, \Dop^\eff := \Dop^\eff(\zbar t), \sigmabf^\eff := \sigmabf^\eff(\zbar t)$; also all other coefficients are evaluated according to the configuration $\zbar \Om$. Therefore, the loads $\zbar\gb$, $\zbar\fb$ and $\zbar\hb$ are given by interpolation
$\zbar\hb:= \hb^{k+1/2}= (\hb^k + \hb^{k+1})$ and in analogy for $\zbar\gb$ and $\zbar\fb$.


The corrector step is based on the approximation
\begin{equation}\label{eq-corr1}
\begin{split}
0 \stackrel{!}{=} \Phi_{k+1}(\Scal(t_{k+1}),\Zcal) \approx \Phi_k(\{\Scal^l\}_{l\leq k},\Zcal) + \delta\zbar\Phi(\zbar{\Scal},\{\Scal^l\}_{l\leq k},\Zcal)\circ (\delta\Scal,\delta \Vcal)\;,
\end{split}
\end{equation}

In the corrector step, we shall use the  approximations according to \eq{eq-cif-31aR} with $t_* = (t_{k+1} +t_k)/2$,  where $\dlt t = t_k-t_{k-1} = t_{k+1} -t_k$,
\begin{equation}\label{eq-corr1a}
\begin{split}
\dot \ub(t_k) & \approx \bar\ub^\pred/\dlt t\;, \quad \dot p(t_k) \approx \bar p^\pred/\dlt t\;,\\
\ddot \ub(t_k) & \approx (\ub^\pred - \bar\ub^\pred)/(\dlt t)^2\;.
\end{split}
\end{equation}

Using \eq{eq-corr1a} the residuals can be introduced
\begin{equation}\label{eq-corr5a}
\begin{split}
\Phi_k^S(\vb) & = \int_{\Om_k}\left (
\sigmabf^k: \nabla \vb - \ol{\rho}^k(\gb^k -(\ub^\pred-\bar\ub^\pred)/(\delta t )^2)\cdot \vb
\right)-  \int_{\Gms_k} \hb^k\cdot\vb\;,\\
\Phi_k^F(q) & = \int_{\Om_k}\left (q\Bb^k:\nabla{\bar\ub^\pred}
+ \delta t \Kb^k(\nabla p^k -\rho_f^k({\fb^{k}} -(\ub^\pred-\bar\ub^\pred)/(\delta t)^2  ))\cdot \nabla q  + q M^k \bar p\right)\;.
\end{split}
\end{equation}


We can now define the corrector problem involving the convection field $\ub^*:= \ub^\pred$ and also the pressure perturbation $p^* := p^\pred$.
Compute the increments $(\ub,p) \in \delta V(t_{k+1}) \times \delta Q(t_{k+1})$ associated with the time step $\dlt t$,  such that
\begin{equation}\label{eq-corr2}
\begin{split}
 &\int_{\Om} \Dop^\eff e(\ub):e(\vb) + \int_{\Om} \sigmabf^\eff :\nabla \vb(\nabla \ub)^T \\
& - \int_{\Om}  \hat p \Bb:\left[(\nabla\vb)(\nabla\cdot \ub)
- \nabla\vb \nabla\ub \right]
- \int_{\Om} (p \Bb +  \hat p \delta \Bb):\nabla\vb \\
& -\int_\Om \left(\zbar\gb -(\ub - \bar\ub)/(\dlt t)^2 )\right)\zbar R(\ub^*,p^*) \cdot \vb
+\int_\Om  \ol{\rho}\ub(\dlt t)^{-2} \cdot \vb \\
& =
\int_\Om  \ol{\rho}\left(\dlt \gb^{k+1} +  (2\bar \ub - \bar\ub^{k-1})/(\dlt t)^2 \right)\cdot\vb + \int_\Gms (\dlt\hb^{k+1} + (\nb\cdot\ub^*)\Hcal \zbar\hb )\cdot\vb -\Phi_k^S(\vb)
\end{split}
\end{equation}
for all $\vb \in V_0(\zbar t)$ and (again multiplied by $\delta t$)
\begin{equation}\label{eq-corr3}
\begin{split}
& \int_{\Om} q \Bb:\left[(\nabla\ub)(\nabla\cdot \ub^*)
- \nabla\ub \nabla\ub^* \right] + \int_{\Om} q \Bb:\nabla\ub
+ \int_{\Om} q M\left( 1 + \nabla\cdot \ub^*\right) p + \int_{\Om} q \delta M^* p
\\
&
+ \delta t\int_{\Om}  \Kb\left(\nabla p -\rho_f\dlt\fb^{k+1} +\rho_f (\ub - 2\bar\ub + \bar\ub^{k-1})/(\dlt t)^2\right) \cdot \nabla q
\\
& -  \delta t\int_{\Om} \Kb \left(\nabla \hat p\nabla \ub^* +\frac{\rho_f}{\kappa_f} p^*(\zbar\fb - (\ub -\bar\ub)/(\dlt t)^2))\right) \cdot \nabla q \\
&+ \delta t\int_{\Om} \Kb \left(\nabla \hat p -\rho_f(\zbar\fb -(\ub - \bar\ub)/(\delta t)^2  )
\right)\cdot (\nabla q\nabla \cdot \ub^*- \nabla q \nabla \ub^*)\\
&+ \delta t\int_{\Om} \delta \Kb (\nabla \hat p -\rho_f(\zbar\fb -(\ub - \bar\ub)/(\delta t)^2  ))\cdot \nabla q \\
& =  \int_{\Om} q\left(\Bb:\nabla \bar\ub - \delta\Bb^*:\nabla\ub^*\right) + \int_{\Om} q M \bar p - \Phi_k^F(q)\;,
\end{split}
\end{equation}
for all $q \in Q_0(\zbar t)$.
Some straightforward rearrangements of these equations lead to the definition of the increment $(\ub^\corr,p^\corr)\equiv(\ub,p)$.

\begin{mydefinition}{def4R}{Weak solution of the corrector subproblem}
The weak solution of the corrector subproblem at time level $t_{k+1}$ is the couple
$(\ub,p) \in \delta V(t_{k+1}) \times \delta Q(t_{k+1})$ which satisfies
\begin{equation}\label{eq-corr6}
\begin{split}
 &
\int_{\Om}(\ol{\rho} + \tilde R(\ub^*,p^*))/(\dlt t)^{2}\ub \cdot \vb
+\int_{\Om}[\tilde \Aop \nabla \ub]:\nabla\vb  - \int_{\Om} p \Bb:\nabla\vb\\
&  =
\int_\Om \left(\tilde R(\ub^*,p^*)(\zbar\gb+ \bar \ub /(\dlt t)^2) + \ol{\rho}(\dlt\gb^{k+1} +  (2\bar \ub - \bar \ub^{k-1}) /(\dlt t)^2 )\right)\cdot\vb  \\
&+\int_\Gms (\dlt\hb^{k+1} + (\nb\cdot\ub^*)\Hcal \zbar\hb )\cdot\vb
+ \int_{\Om}\hat p \delta\Bb :\nabla\vb -\Phi_k^S(\vb)
\;,
\end{split}
\end{equation}
for all $\vb \in V_0(t_k)$ and
\begin{equation}\label{eq-corr7}
\begin{split}
& \int_{\Om} q [\Bb + \tilde\Bb(\ub^*)]:\nabla\ub
 +  \delta t\int_{\Om} \Kb \nabla p \cdot \nabla q
 + \int_{\Om} q (M + \tilde M(\ub^*) + \delta M^*) p\\
&+ (\dlt t)^{-1}\int_{\Om} \rho_f\nabla q \cdot \left(\Kb + \what\Kb(\ub^*,p^*) + \dlt\Kb\right)\ub
\\
& =   \int_{\Om} q\left(\Bb:\nabla\bar\ub - \delta\Bb :\nabla \bar\ub + M\bar p\right)  
- \delta t \int_{\Om}
\left (\tilde\Kb(\ub^*) +  \delta \Kb\right)\nabla \hat p \cdot \nabla q\\
&
+ \delta t \int_{\Om}\rho_f\nabla q \cdot[\what\Kb(\ub^*,p^*) + \delta \Kb]\rb^A + \delta t \int_{\Om}\rho_f\nabla q \cdot\Kb \rb^B
 - \Phi_k^F(q)
\;,
\end{split}
\end{equation}
for all $q \in Q_0(t_k)$,
where 
\begin{equation}\label{eq-corr8}
\begin{split}
\zbar \rb^A  = \zbar \fb + \bar\ub /(\delta t)^2 \;,\quad
\zbar \rb^B  = \dlt\fb^{k+1} + (2\bar\ub - \bar\ub^{k-1})/(\dlt t)^2\;.
\end{split}
\end{equation}

\end{mydefinition}

The new configuration $\Om^{k+1}$ is established
using the increments $\ub^\corr := \ub$ and $p^\corr:= p$ associated with the time step $\dlt t$.
Consequently, although the mid-time level states are used only in the corrector problem, they can be updated by an interpolation such as the following one which has been employed in this work:
\begin{equation}\label{eq-corr8b}
\begin{split}
\ub^{k+1/2} = \ub^k + (\ub^\pred + \ub^\corr)/2\;,\quad
p^{k+1/2} = p^k + (p^\pred + p^\corr)/2\;.
\end{split}
\end{equation}

When computing increments $(\ub^1,p^1)$, at time level $k=0$, $\bar\ub := \delta t \dot\ub(0)$ can be defined due to the initial condition, whereas $\bar\ub^{-1}$ is undefined. However, $\bar\ub^{-1}$ occurs in \eqs{eq-corr6}{eq-corr8} involving the expression
$\bar\ub - \bar\ub^{-1} \approx  (\delta t)^2\ddot\ub(0)$. It approximates the initial acceleration which can be computed easily due to the linearity of the residual $\Phi_t$  at time $t=0$  \wrt acceleration  and pressure rate fields $\ddot\ub(0)$ and $\dot p(0)$. For this we consider \eq{eq-cif-6R} with the stress defined according to \eq{eq:1} where the effective stress is expressed using the given initial conditions on the displacement, $\ub^0 = \ub(0)$, and pressure $p^0 = p(0)$. At time $t=0+$, the model equations in \eq{eq-cif-8} read as
\begin{equation}\label{eq-cif-00-BR}
\begin{split}
\int_\Om \ol{\rho}\ddot\ub^{0+}\cdot\vb & = \int_{\pd\Om}\hb^{0+}\cdot\vb + \int_\Om \ol{\rho}\gb^{0+}\cdot\vb - \int_\Om\sigmabf(\ub^0,p^0):\nabla\vb \quad \forall\vb\in V_0(0)\;,\\
\int_\Om q M \dot p^{0+} & = - \int_\Om\nabla q \cdot \Kb\left(\nabla p^0 + \rho_f(\ddot\ub^{0+} - \fb^{0+})\right) - \int_\Om q \Bb:\nabla\dot\ub^0\quad \forall q \in Q_0(0)\;,
\end{split}
\end{equation}
where $\ub(0),\dot\ub(0)$ and $p(0)$ are given by the initial conditions. First $\ddot\ub^{0+}$ can be solved using \eq{eq-cif-00-BR}$_1$, consequently $\dot p^{0+}$ is obtained from
\eq{eq-cif-00-BR}$_2$. Then we put $\bar\ub^{-1} = \bar\ub - \ddot\ub^{0+}(\delta t)^2$ and also $\bar p = \dot p^{0+}\delta t$.

\subsection{Spatial discretization and computational algorithm}
Problem \eqs{eq-pred6}{eq-pred7} and \eqs{eq-corr6}{eq-corr7} can be solved using the finite element method (FEM).
Since details related to the numerical approximation are beyond the scope of the present publication, we focus on the solution algorithm. For the sake of clarity we shall return to the labeling the variables associated with time $t_k$ by the superscript $^k$ and denote the finite time increments by $\dlt$, thus, $\ub^{k+1} = \ub^k + \dlt \ub$.

For the spatial discretization we use the mixed conforming finite elements with piecewise $Q2$ approximation for the displacements and $Q1$ approximation for the fluid pressure; this combination which is coherent with the restriction of the Babu\v{s}ka--Brezzi condition \cite{BrezziBook1991} is quite usual when dealing with this kind of mixed displacement--pressure formulations. For both the incremental subproblems, \ie the predictor and the corrector problems, the discretization yields the following linear equations (a self-explaining notation is used):

\begin{equation}\label{eq-num-1}
\begin{split}
[(\dlt t)^{-2} \Mbm + \Abm] \dlt\ubm  - \Bbm \dlt\pbm & = \fbm\;,\\
[(\dlt t)^{-1} \Nbm + \Cbm] \dlt\ubm  + [\dlt t \Kbm   +  \Dbm] \dlt\pbm & = \gbm\;,
\end{split}
\end{equation}
where all the matrices and the \rhs column matrices are constituted by tensors and vectors
given at time $t_k$. For the assumed boundary conditions matrices $\Abm$, $\Kbm$, $\Dbm$ and $\Mbm$ are symmetric positive definite. In the linear case, \eq{eq-num-1} is obtained from the time discretization of the Biot model, so that $\Cbm = \Bbm^T$.
In the non-linear case, while  $\Bbm$ is unchanged, the matrices $\Cbm$ arising from the first integral in \eq{eq-pred7}, or \eq{eq-corr7} reflect the modified tensor $\Bb + \tilde\Bb(\ub^*)$, thus $\Cbm$ is only similar to $\Bbm^T$.
If an incompressible medium is considered, thus both the solid and the fluid phases are incompressible, the coefficient $M$ vanishes, so that $\Dbm = \textbf{0}$. Therefore, to allow for a small step $\dlt t$, a stable approximation satisfying the ``inf-sup'' condition associated with the mixed formulations known from the elasticity should be used, see \eg \cite{BrezziBook1991}.

We recall the simplifying assumption of the zero initial conditions, including the first order time derivatives, so that $\ub,\dot\ub,p,\dot p$ vanish at $t_0$.
Below by $\Om_k$ we mean the domain partitioned by the FEM at time $t_k$; thus, $\Om_k$ can be considered as the FEM mesh.
The flowchart of the computational algorithm consists of the following steps.

\begin{enumerate}
\item Initialization: put $k=0$, in $\Om_0$, set the initial conditions $\ub^0,p^0$ and also put $\dlt\ub^0 = \dot\ub(0)\dlt t$.
By solving \eq{eq-cif-00-BR}, compute $\dlt \ub^{-1} = \dlt\ub^0 - \ddot\ub^{0+}(\dlt t)^2$ and $\dlt p^0 = \dot p^{0+}\delta t$.

\item Update the reference configuration $\Om_k$ at $t_k$ and define the incremental model  with load increments $\dlt \gb^{k+1},\dlt \gb^{k+1}$, and $ \dlt \hb^{k+1}$ defined at FE mesh nodes.

\item Compute $(\ubm^\pred,\pbm^\pred)$ by solving the discretized Predictor problem \eqs{eq-pred6}{eq-pred7} with load increments $\dlt \gb^{k+1},\dlt \gb^{k+1}$, and
$ \dlt \hb^{k+1}$

\item Introduce the intermediate configuration $\zbar \Om$ at $t_{k+1/2}$ where $\ubm^{k+1} = \ubm^k + \ubm^\pred/2$ and
$\pbm^{k+1} = \pbm^k + \pbm^\pred/2$. Update the FE mesh.

\item Compute $(\ubm^\corr,\pbm^\corr)$ by solving the discretized Corrector problem \eqs{eq-corr6}{eq-corr7} with updated load increments $\dlt \gb^{k+1},\dlt \gb^{k+1}$, and $ \dlt \hb^{k+1}$ defined at FE mesh nodes.

\item Update the solution: $\ubm^{k+1} = \ubm^k + \ubm^\corr/2$ and $\pbm^{k+1} = \pbm^k + \pbm^\corr/2$.
Correct the intermediate configuration using \eq{eq-corr8b}.

\item If $k < \bar k$, put $k:=k+1$ and go to step 2.
\item Final postprocessing.

\end{enumerate}

The procedure executed in steps 2 and 4 of the global algorithm comprise the following actions:
\begin{itemize}
\item Given the total displacement and pressure fields in terms of the nodal values $\ubm$ and $\pbm$, respectively, and the convection field $\ubm^*$ at each FE integration point of $\Om$ do:
\begin{itemize}
\item Recover the gradients $\nabla \ub, \nabla \ub^*, \nabla p$ in $\Om$ and $\nabla_X \ub$ in $\Om_0$.
\item Constitute the deformation gradient $\Fb = \Ib + \nabla_X \ub$, and $J = \det \Fb$.
\item Compute the effective stress $\sigmabf^\eff$ using \eq{eq:1} update the model incremental parameters $\Dop^\eff,\Bb,\Kb,M$.
\item Compute the variations $\dlt \Kb, \dlt \Bb$ and $\dlt M$ associated with $\ub^*$ and $\nabla  \ub^*$.
\end{itemize}
\item Assemble the matrices involved in \eq{eq-num-1} of the FEM incremental model (for the predictor, the corrector)
\end{itemize}

The perturbation $\dlt\Bb = \bmi{0}$ due to the solid incompressibility assumed in this paper. To compute $\dlt \Kb$ and
$\dlt M$, the convection increment $\dlt\ub^*$ is employed in \eq{eq-cif-29d-BR}.


\section{Numerical examples}\label{sec:ex}

\def\figdir{}

In this section we present 2D examples which illustrate the model performance.
The numerical solution is based on the mixed FEM with conforming Q2-Q1 approximation.
The model has been implemented within our in-house developed code SfePy~\cite{cimrman_2014:sfepy}

\subsection{Model validation}

In order to validate the proposed non-linear model of the FSPM, we follow two simulation tests
performed in Li et al. \cite{li_borja_regueiro_2004}. For this purpose we also implemented the
infinitesimal formulation of the dynamic problem of porous media including
the Newmark integration scheme, as the one used in the referenced article.

\subsubsection{Test \rom{1}: compression}

\begin{figure}
  \hfil\includegraphics[width=0.9\linewidth]{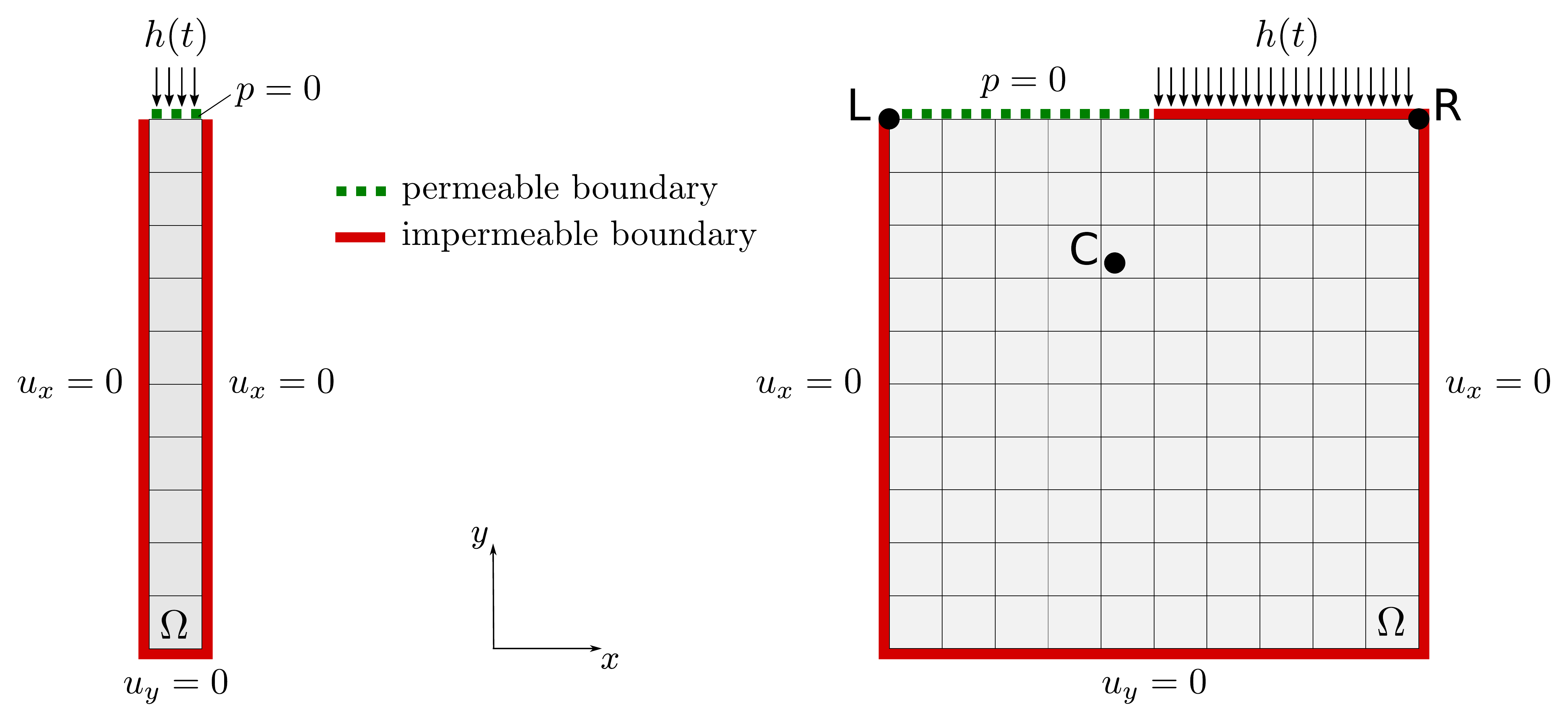}\hfil
  \caption{Finite element meshes and applied boundary conditions
    used for model validation. The green dotted line represents
    the permeable part of the boundary while the rest of the surface
    is assumed to be impermeable (red solid line).
    The roller supports are applied to the left, right and bottom
    boundaries, $h(t)$ is a surface load.}
    \label{fig_val1}
\end{figure}

In the first test, we compare responses of the following three models: 1) the linear model solved by the Newmark time
integration scheme, cf. \cite{li_borja_regueiro_2004}, 2) the linear model where the time derivatives are approximated
using the backward finite difference method, and 3) the non-linear model with the predictor-corrector time integration
described in the previous sections.  The hyperelastic model \eq{eq:1b} of the
solid skeleton provides also the Hookean material model for the linear problems.
 The rectangular domain 1$\times$10\,m is partitioned using a finite element mesh
consisting of 10 elements of dimensions 1 $\times$ 1\,m is shown in
Fig~{\ref{fig_val1}}~left.
In what follows, by lower indices $x,y$ we refer to the coordinate axes.
 We apply the following boundary conditions:
$u_x = 0$ at left and right edges, $u_y = 0$ at the bottom edge and $p = 0$ at the
top edge, so the body surface except the top part is assumed to be
impermeable. The top boundary is loaded by boundary traction $h_y(t) = - H(t) \bar h$,
where $H(t)$ is the Heaviside step function and $\bar h = 40$\,kPa.

\begin{figure}
  \hfil
  \includegraphics[width=0.49\linewidth]{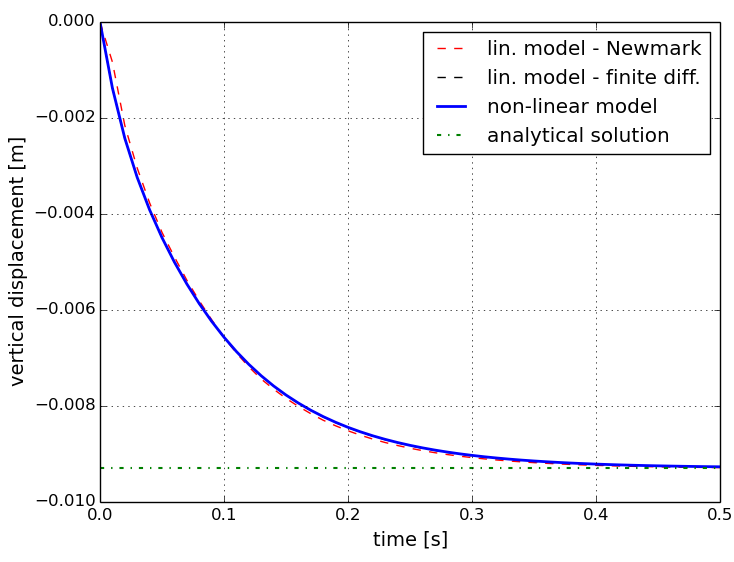}
  \includegraphics[width=0.49\linewidth]{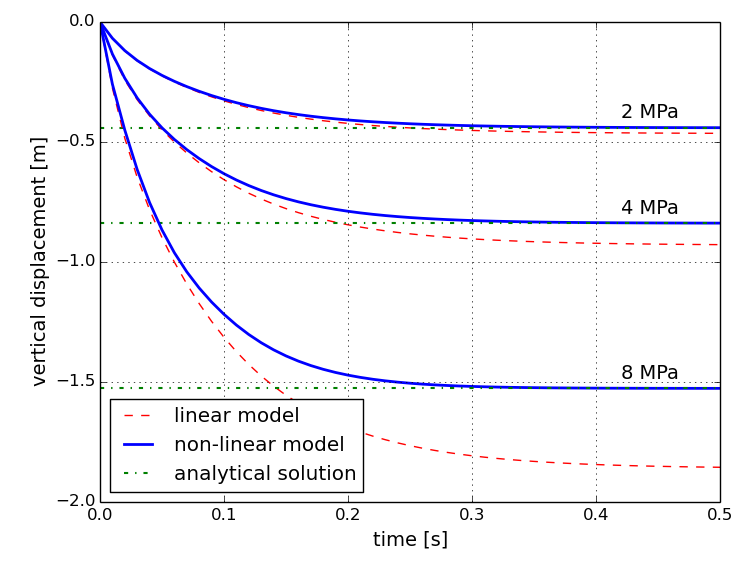}
  \hfil
  \caption{Compression test -- time evolution of the vertical displacement: left)
   load $\bar h = 40$\,kPa, right) $\bar h = 2, 4, 8$\,MPa.}
    \label{fig_val2}
\end{figure}

The resulting time evolution of the vertical displacement of the top surface is
shown in Fig.~{\ref{fig_val2}}~left. The curves corresponding to the linear and
non-linear models solved by the finite differences coincide and are almost
identical to the curve representing the response of the linear model solved by
the Newmark scheme. The time step in all cases was $\delta t = 0.01$\,s. The
analytic solution of the steady-state compression is $u_y = e_y H_m$, where
$H_m$ is the mesh height and the vertical strain $e_y$ can be obtained by
solving the non-linear equation:
$$
(1 + e_y)\mu + {1\over 1 + e_y}\left(\lambda \ln(1 + e_y) - \mu\right) + \bar h = 0.
$$

The material parameters of the solid skeleton are as follows: Lam\'e
coefficients $\lambda = 29$\,MPa, $\mu = 7$\,MPa; mass density $\rho_{s0} =
2700$\,kg/m$^3$ and initial volume fraction $\Phi^{s0} = 0.58$. The fluid part
is characterized by: mass density $\rho_{f0} = 1000$\,kg/m$^3$, bulk stiffness
$K_f = 22$\,GPa and hydraulic permeability $\Kb = K_0 \Ib e^{\beta(J - 1)}$,
where $K_0 = 0.1 / (g \rho_{f0})$\,m$^{2}$ / (Pa $\cdot$ s) and $\beta = 0.8$.
In Fig.~{\ref{fig_val2}}~right, we compare the non-linear solutions to the
corresponding linear solutions and to the analytical predictions of non-linear
response for surface loads $\bar h = 2, 4, 8$\,MPa.

\subsubsection{Test \rom{2}: partial compression}\label{sec_test2}

\begin{figure}
  \hfil
  \includegraphics[width=0.6\linewidth]{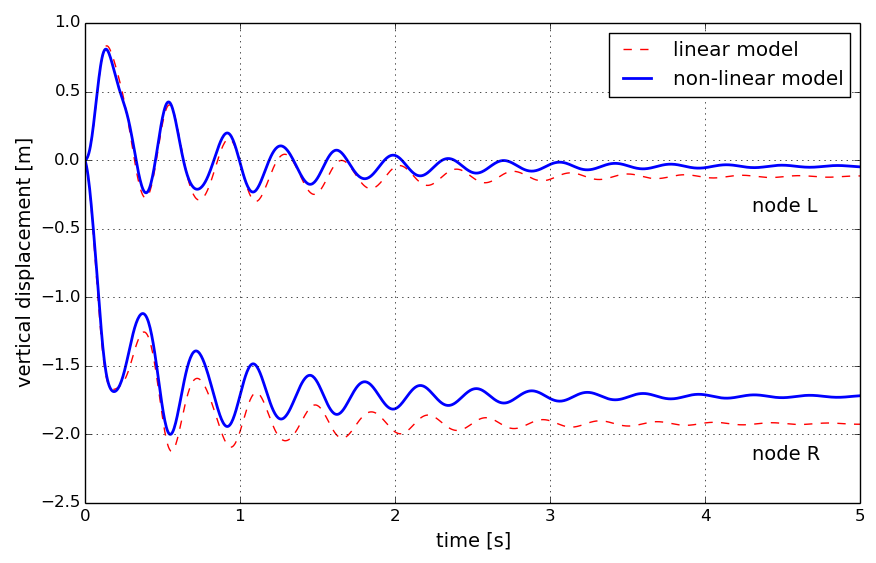}
  \hfil
  \caption{Partial compression test -- time evolution of the vertical
   displacement at nodes L and R, see Fig.~{\ref{fig_val1}}~right.}
    \label{fig_val3}
\end{figure}

Now we consider the square-shaped domain $10\times10$\,m  partitioned by
$10\times10$ elements. The boundary conditions, as shown in
Fig.~{\ref{fig_val1}}~right, are similar as in the first test, however, the upper boundary is loaded on its right half only which is also non-permeable. The fluid can evacuate through its left part, where $p = 0$ is prescribed. The material parameters are the same as in the
previous test except the Lam\'e coefficients which are now $\lambda = 8.4$\,MPa,
$\mu = 5.6$\,MPa and the time step is $\delta t = 0.005$\,s.
Fig.~{\ref{fig_val3}} shows the time responses of the linear and non-linear
models under partial compression $\bar h = 4$\,MPa at two mesh nodes L and R,
see Fig.~{\ref{fig_val1}}~right. The non-linear response is smaller than the
linear one which is in agreement with the previous results. When decreasing the
applied load $\bar h$ the both responses coalesce.

\subsection{Effect of time step}

\begin{figure}
  \hfil
  \includegraphics[width=0.99\linewidth]{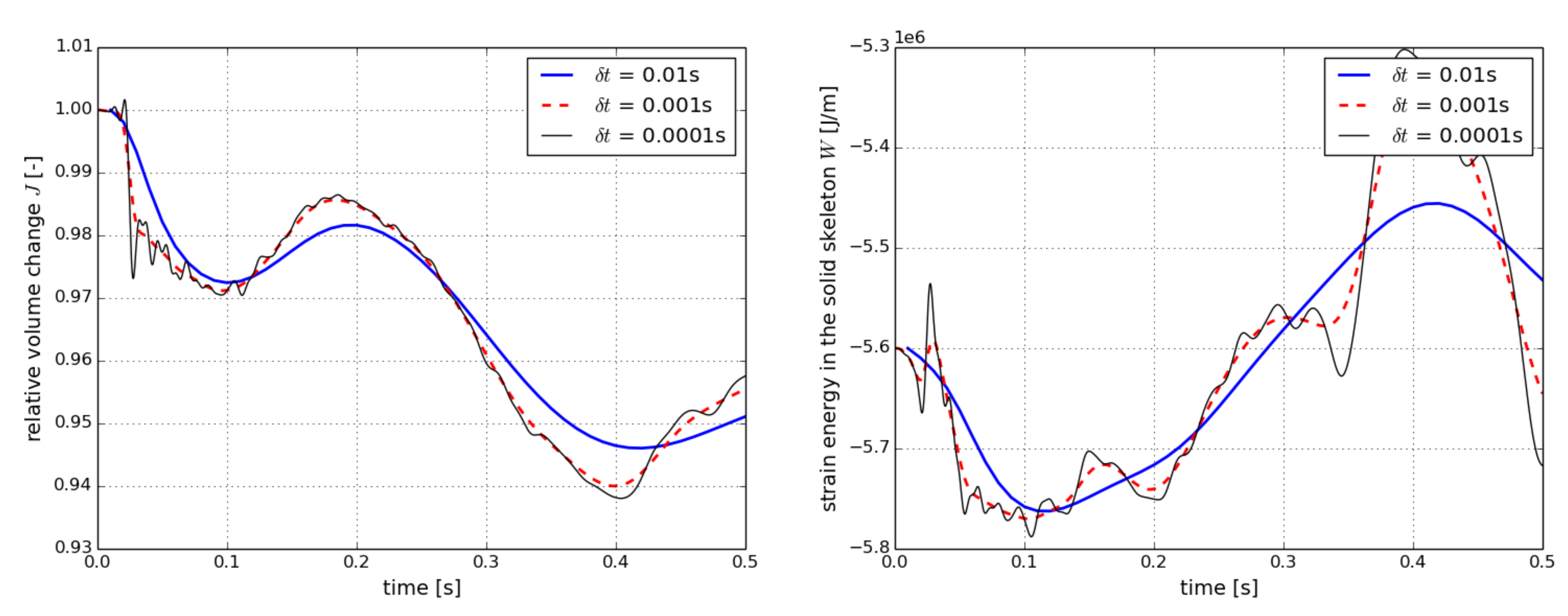}
  \hfil
  \caption{
    Relative volume change $J$ and hyperelastic strain energy $W$
    computed for time increments $\delta t = 0.01, 0.001, 0.0001$\,s
    at point C of domain $\Omega$.}
    \label{fig_simul1}
\end{figure}

\begin{figure}
  \hfil
  \includegraphics[width=0.47\linewidth]{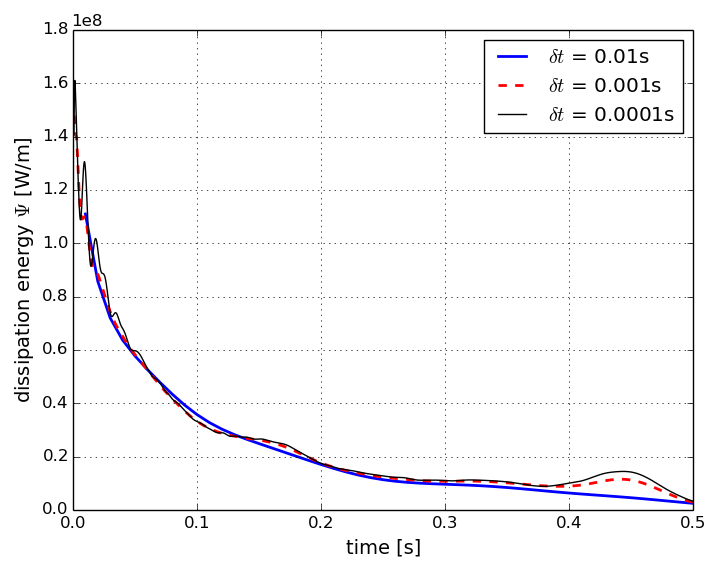}
  \hfil
  \caption{
    Dissipation energy $\Psi$ computed for time increments
    $\delta t = 0.01, 0.001, 0.0001$\,s.}
    \label{fig_simul3}
\end{figure}

\begin{figure}
  \hfil
  \includegraphics[width=0.99\linewidth]{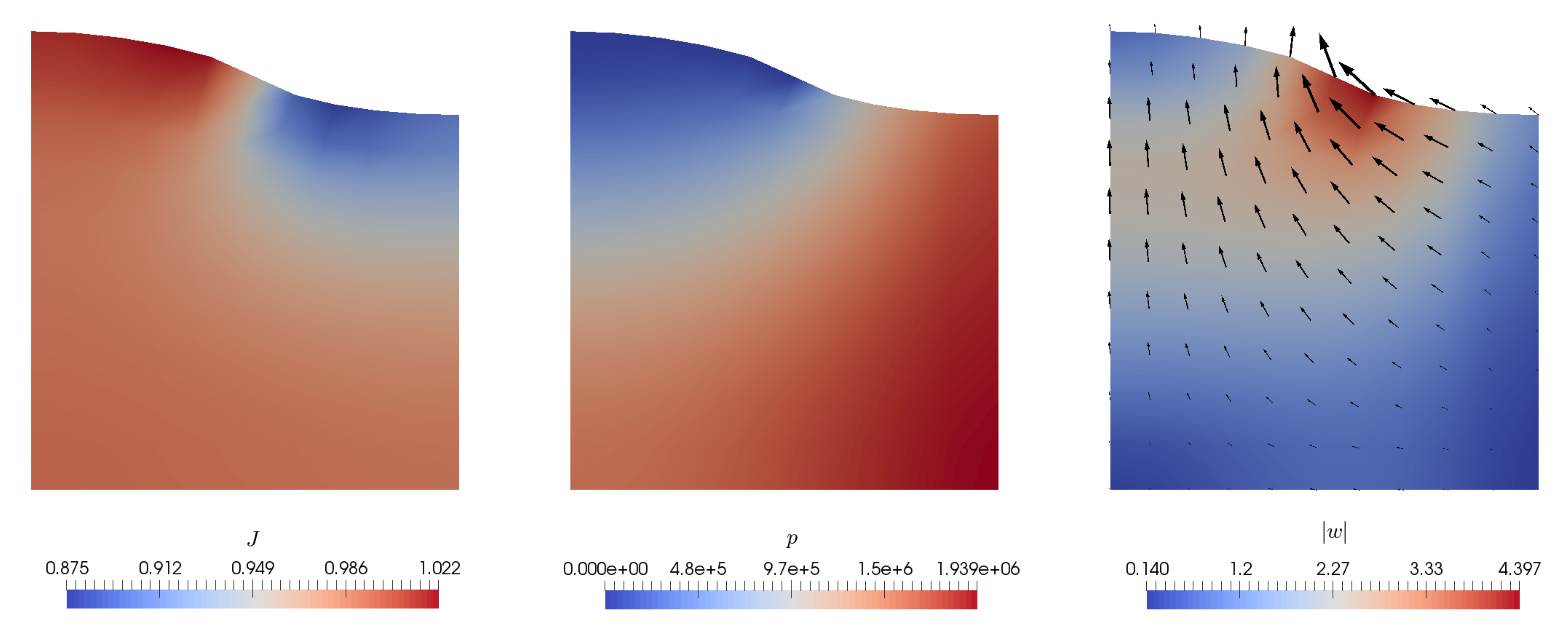}
  \hfil
  \caption{
    Distribution of relative volume change $J$ (left),
    fluid pore pressure $p$ (middle)
    and seepage velocity magnitude $\vert\bmi{w}\vert$ (right)
    at time $t = 0.1$\,s ($\delta t = 0.01$\,s).}
    \label{fig_simul2}
\end{figure}

In contrast with nonlinear elasticity problems, where the equilibrium can be
ensured by reiterations while updating the reference configuration, the choice
of a convenient time step in problems featured by the dissipation is an
important issue, since such reiterations are not easily justified. On one hand,
the smaller the step, the smaller the error induced by all the approximations
associated with the time integration. On the other hand, the time step is
also bounded from below to prevent the stability. In this part we solve the problem
identical to the one defined in Sec.~\ref{sec_test2} including the boundary
conditions and material parameters. The dependency of the solution on $\delta t$
is illustrated for three different time steps in Figures~\ref{fig_simul1} and
\ref{fig_simul3}.
In Fig.~\ref{fig_simul1}, we compare the relative volume change $J$ and the
hyperelastic strain energy $W$ computed for time increments $\delta t = 0.01,
0.001$ and $0.0001$\,s, the values at point C of domain $\Omega$, see
Fig.~{\ref{fig_val1}}~right, are depicted. The influence of the time step length
on the dissipation energy $\Psi = \int_\Omega \Kb \nabla p \cdot \nabla p$ is
shown in Fig.~\ref{fig_simul3}. In general, when decreasing the time step $\delta t$ without any FE mesh refinement of the spatial discretization, numerical
oscillations may deteriorate the solution, since for $\dlt t \rightarrow 0$ the numerical solution converges to the exact solution of the semidiscretized problem in space.
 This well known effect pointed out \eg in \cite{li_borja_regueiro_2004} is apparent
from Fig.~\ref{fig_simul1} where remarkable oscillations of computed responses
obtained with time increment $\delta t = 0.0001$\,s appear. The relative volume change,
pore pressure field and seepage velocity magnitude at time $t = 0.1$\,s ($\delta
t = 0.01$\,s) are presented in Fig.~\ref{fig_simul2} where also the deformed FE
mesh is shown.

\section{Conclusion}\label{sec:concl}
In this work we derived an incremental formulation for computational analysis of the fluid-saturated porous media undergoing large deformation. As a modelling assumption, we consider the Biot model describing the fluid-structure interaction with respect to perturbations in the spatial Eulerian configuration at the macroscopic level. The rate form of the equilibrium and the mass conservation equations is derived. This leads to the incremental formulation coherent with the updated Lagrangian formulation.
The time discretization and approximations of the time derivatives using the increments lead to
the predictor--corrector algorithm where the latter uses a secant approximation of the residual function differential.
Further modifications of this algorithm are possible towards higher order implicit integration schemes introducing intermediate steps (multi-stage algorithms), cf. \cite{Carstens-Kuhl-time-integr-AAM2012}
The relationship between the FEM spatial approximation and the time step associated with the CFL condition is an important issue which requires to estimate the speed of the wave propagation.

The proposed nonlinear Biot-type model has been implemented for 2D problems in the SfePy code \cite{cimrman_2014:sfepy}. A thorough validation of the model and verification of the incremental formulation is a complex task which exceeds scopes of this publication, nevertheless, we have shown that the proposed computational model leads to results which are in accordance with the published work \cite{li_borja_regueiro_2004}.

As mentioned in the introduction, this work is intended as a basis for our further research which will focus on the multiscale modelling of the porous material subject to large deformations. Since the direct homogenization of the nonlinear problem is cumbersome,
it is advisable to benefit from a suitable combination of both the phenomenological and upscaling approaches which
provide a suitable tool for a computational homogenization. In particular,  the homogenization of the fluid-structure interaction can provide the effective medium poroelastic coefficients and the permeability as well as the the sensitivity of these material parameters involved in the incremental formulation. In a series of papers \cite{Rohan2003,AMC2015-Rohan-Lukes} and \cite{Rohan-etal-CMAT2015} we tackled this topic, however the problem is still open.

\paragraph{Acknowledgments:} This research is supported by the project
 LO 1506 of the Czech Ministry of Education, Youth and Sports.

\end{document}